\documentclass[11pt]{amsart}
\usepackage{amsmath, amssymb}
\usepackage{amsfonts}
\usepackage[arrow,matrix,curve,cmtip,ps]{xy}

\usepackage{amsthm}

\allowdisplaybreaks

\newtheorem{theorem}{Theorem}[section]
\newtheorem{lemma}[theorem]{Lemma}
\newtheorem{proposition}[theorem]{Proposition}
\newtheorem{corollary}[theorem]{Corollary}

\newtheorem*{theorem*}{Theorem}
\theoremstyle{remark}
\newtheorem{remark}[theorem]{Remark}
\newtheorem{definition}[theorem]{Definition}
\newtheorem{example}[theorem]{Example}

\newtheorem*{acknowledgment}{Acknowledgment}

\numberwithin{equation}{section}

\newcommand{\Z}{\mathbb{Z}}
\newcommand{\N}{\mathbb{N}}

\newcommand{\T}{\mathbb{T}}
\newcommand{\im}{\operatorname{im }}

\newcommand{\aut}{\operatorname{Aut}}

\newcommand{\dirlim}{\operatorname{\displaystyle \lim_{\longrightarrow}}}
\newcommand{\Li}{\mathcal{L}}
\newcommand{\K}{\mathcal{K}}
\newcommand{\M}{\mathcal{M}}

\newcommand{\clspan}{\operatorname{\overline{\mathrm{span}}}}

\newcommand{\TX}{\mathcal{T}_X}
\newcommand{\OX}{\mathcal{O}_X}

\begin{document}
\title{Strong Shift Equivalence of $C^*$-correspondences}

\author{Paul S. Muhly}

\author{David Pask}

\author{Mark Tomforde}

\address{Department of Mathematics\\ University of Iowa\\
Iowa City\\ IA 52242-1419\\ USA}

\email{pmuhly@math.uiowa.edu}

\address{School of Mathematical and Physical Sciences\\ The University of Newcastle\\ University Drive, Callaghan\\ NSW 2308, AUSTRALIA}

\email{david.pask@newcastle.edu.au}

\address{Department of Mathematics\\ The College of William and Mary\\ P.O. Box 8795\\ Williamsburg, VA 23187-8795\\ USA}

\email{tomforde@math.wm.edu}

\thanks{The first author was supported by NSF Grant DMS-0355443 and the
third author was supported by NSF Postdoctoral Fellowship DMS-0201960.}

\date{\today}


\begin{abstract}
We define a notion of strong shift equivalence for $C^*$-correspondences and show that strong shift equivalent $C^*$-correspondences have strongly Morita equivalent Cuntz-Pimsner algebras. Our analysis extends the fact that strong shift equivalent square matrices with non-negative integer entries give stably isomorphic Cuntz-Krieger algebras.

\end{abstract}

\maketitle

\section{Introduction}

Inspiration for this work comes in large part from the two papers, \cite{rW73} by R. Williams and \cite{CK80} by J. Cuntz and W. Krieger.  To understand why, suppose $A$ and $B$ are two square matrices, possibly of different sizes, whose entries are non-negative integers.  Then $A$ and $B$ are called \emph{strong shift equivalent} if there is a finite chain of square matrices with non-negative integer entries, $A_1, A_2, \ldots A_n$, such that $A_1 = A$, $A_n = B$ and for each $i$, $1 \leq i \leq n-1$, there is a pair of matrices $(R_i, S_i)$ with non-negative integer entries such that $A_i = R_iS_i$ and $A_{i+1} = S_iR_i$. (An individual pair $(R_i,S_i)$ is sometimes called an elementary strong shift equivalence between $A_i$ and $A_{i+1}$.) On the other hand, if $A$ is a square matrix with non-negative integer entries, then there is a well-known process for building a shift dynamical system $(X_A, \sigma_A)$ --- a so-called shift of finite type: One regards $A$ as the incidence matrix of a finite graph, $E = (E^0,E^1,r,s)$, where $E^0$ is the space of vertices, $E^1$ is the space of edges, and $r$ (resp.~$s$) is the range map (resp.~source map) from $E^1$ to $E^0$. The space $X_A$, then, is the two-sided infinite path space $\{(e_i)_{i\in \mathbb{Z}} \in ({E^1})^{\mathbb{Z}} \mid s(e_{i+1}) = r(e_i) \}$. Evidently, $X_A$ is a closed subset of the compact space $({E^1})^{\mathbb{Z}}$ that is invariant under the shift map $\sigma$, given by $\sigma (e_i) = (e_{i+1})$. The shift map $\sigma$ is a homeomorphism, of course, and therefore, so is its restriction to $X_A$, which is denoted $\sigma_A$. In \cite[Theorems~A \& F]{rW73} Williams proved that the shift dynamical systems $(X_A,\sigma_A)$ and $(X_B,\sigma_B)$ are conjugate, meaning that there is a homeomorphism $\phi$ from $X_A$ to $X_B$ such that $\phi \circ \sigma_A = \sigma_B \circ \phi$, if and only if $A$ and $B$ are strong shift equivalent.  

Subsequently, in \cite{CK80}, Cuntz and Krieger attached a $C^*$-algebra, $\mathcal{O}_A$, to every square matrix $A$ having non-negative integer entries. Cuntz and Krieger worked primarily with matrices whose entries have zeros and ones, but in \cite[Remark~2.16]{CK80} they observe that their construction for zero-one matrices can be adjusted to cover matrices with non-negative integer entries.  In \cite[Theorem~3.8]{CK80} they show that if $(X_A,\sigma_A)$ and $(X_B,\sigma_B)$ are conjugate, then $\mathcal{O}_A$ and $\mathcal{O}_B$ are strongly Morita equivalent in the sense of Rieffel.  Coupled with Williams's theorem, we conclude that if $A$ and $B$ are two strong shift equivalent matrices, then the $C^*$-algebras $\mathcal{O}_A$ and $\mathcal{O}_B$ are strongly Morita equivalent.

Reflecting on the graphs associated with shift dynamical systems and taking into account developments in the theory of Cuntz-Krieger algebras that allow one to express an $\mathcal{O}_A$ in terms of a graph (see \cite{Rae}, for example), we were led to consider the following generalization of strong shift equivalence and its relation to so-called Cuntz-Pimsner algebras.  (Precise definitions of the terms used and technical hypotheses, which are omitted here, will be developed thoroughly in the body of the paper.)  Suppose $A$ and $B$ are $C^*$-algebras and that $E$ (resp.~$F$) is a $C^*$-correspondence over $A$ (resp.~$B$) (so in particular, $E$ and $F$ are bimodules over $A$ and $B$, respectively).  Then we shall say that $E$ is (elementary) strong shift equivalent to $F$ in case there is a correspondence $R$ from $A$ to $B$ and a correspondence $S$ from $B$ to $A$ (in particular, $R$ is an $A$--$B$-bimodule and $S$ is a $B$--$A$-bimodule)  such that $E \cong R \otimes_B S$ and $F \cong S \otimes_A R$. Given the relation between strong shift equivalence of matrices and the strong Morita equivalence of the associated Cuntz-Krieger algebras, we were led to speculate that if $E$ and $F$ are strong shift equivalent, then the Cuntz-Pimsner algebras $\mathcal{O}_E$ and $\mathcal{O}_F$ are strongly Morita equivalent.  It turns out that our speculation is correct, at least under appropriate hypotheses, as we shall show in Theorem \ref{SSE-correspondences}.  Our result captures the connection just discussed between strong shift equivalence of matrices and the strong Morita equivalence of their associated Cuntz-Krieger algebras, as we shall make clear in Section \ref{Examples and Counterexamples}.

In the next section we develop the basic facts about correspondences, Cuntz-Pimsner algebras, etc. that we need.  We note that most places in the literature, when constructing Cuntz-Pimsner algebras of correspondences, the blanket assumption is made that the coefficient algebra acts faithfully on the correspondence.  However, in our investigation it is important to allow non-faithful actions.  Fortunately, Katsura recently has developed the theory of $C^*$-correspondences where the action need not be faithful in \cite{Kat3}, \cite{Kat4}, and \cite{Kat5} and he has extended the notion of Cuntz-Pimsner algebra to this setting. The next section expands upon his work for our purposes.

In Section \ref{ESSERCC} we prove our main result, Theorem \ref{SSE-correspondences}, just described.   In Sections \ref{non-ess-sec} and \ref{Examples and Counterexamples} we explore the limitations and necessity of some of our technical hypotheses. Also, in Section \ref{Examples and Counterexamples} we show how our analysis relates to graph $C^*$-algebras and shifts of finite type.

\begin{acknowledgment}
We are very grateful to Baruch Solel for valuable insights and conversations that helped with our investigation.  We are also grateful to Berndt Brenken for some very helpful corrections to an initial draft.
\end{acknowledgment}

\section{Preliminaries}\label{preliminaries}

\subsection{$C^*$-correspondences}

We follow the conventions of Lance \cite{Lan} for our terminology of Hilbert $C^*$-modules, and we use the notation and conventions of Katsura in \cite{Kat3}, \cite{Kat4}, and \cite{Kat5} for $C^*$-correspondences.

\begin{definition}
If $A$ is a $C^*$-algebra, then a \emph{right Hilbert $A$-module} is a Banach
space $X$ together with a right action of $A$ on $X$ and an $A$-valued inner
product $\langle \cdot , \cdot \rangle_X$ satisfying
\begin{enumerate}
\item[(i)] $\langle \xi, \eta a \rangle_X =  \langle \xi, \eta \rangle_X a$
\item[(ii)] $\langle \xi, \eta \rangle_X =  \langle \eta, \xi \rangle_X^*$
\item[(iii)]  $\langle \xi, \xi \rangle_X \geq 0$ and $\| \xi \| = \| \langle
\xi, \xi \rangle_X \|^{1/2}$
\end{enumerate}
for all $\xi, \eta \in X$ and $a \in A$.
For a Hilbert $A$-module $X$ we let $\Li(X)$ denote the
$C^*$-algebra of adjointable operators on $X$, and we let $\K (X)$ denote
the closed two-sided ideal of compact operators given by $$\K (X) := \clspan
\{ \Theta_{\xi,\eta}^X : \xi, \eta \in X \}$$ where $\Theta_{\xi,\eta}^X$ is
defined by $\Theta_{\xi,\eta}^X (\zeta) := \xi \langle \eta, \zeta
\rangle_A$.  When no confusion arises we shall often omit the superscript and
write $\Theta_{\xi,\eta}$ in place of $\Theta_{\xi,\eta}^X$.
\end{definition}

\begin{remark} \label{direct-sum-right}
If $X$ is a right Hilbert $A$-module and $Y$ is a right Hilbert $B$-module, then we may give $X \oplus Y$ the structure of a right Hilbert $A \oplus B$-module by defining $(x,y) (a,b) := (xa,yb)$ and $\langle (x_1,y_1), (x_2,y_2) \rangle_{X \oplus Y} := (\langle x_1,x_2\rangle_X, \langle y_1, y_2 \rangle_Y)$.
\end{remark}

\begin{definition}
If $A$ and $B$ are $C^*$-algebras, then a \emph{$C^*$-correspondence from $A$ to $B$} is a right
Hilbert $B$-module $X$ together with a $*$-homomorphism $\phi_X : A \to
\Li(X)$.  We consider $\phi_X$ as giving a left action of $A$ on $X$ by setting
$a \cdot x := \phi_X(a) x$.  When $X$ is a $C^*$-correspondence from $A$ to $B$ we will sometimes write ${}_AX_B$ to keep track of the $C^*$-algebras.  If $A=B$ we refer to $X$ as a \emph{$C^*$-correspondence over $A$}.
\end{definition}

\begin{definition}
Let $X$ and $Y$ be $C^*$-correspondences from $A$ to $B$.  An \emph{isomorphism} from $X$ to $Y$ is a surjective linear map $T : X \to Y$ with the property that $T(x b) = T(x) b$, $\langle T(x), T(y) \rangle_Y = \langle x, y \rangle_X$, and $T(\phi_X(a) x) = \phi_Y(a) T(x)$ for all $x,y \in X$, $b \in B$, and $a \in A$.  We say that $X$ and $Y$ are isomorphic if there exists an isomorphism from $X$ to $Y$, and in this case we write $X \cong Y$.
\end{definition}

Evidently, since $\langle T(x), T(y) \rangle_Y = \langle x, y \rangle_X$, an isomorphism is automatically injective.  Thus isomorphisms are bijective.

\begin{lemma} \label{iso-phi-translate}
Let ${}_AX_B$ and ${}_AY_B$ be $C^*$-correspondences from $A$ to $B$.  Then $T$ induces a $*$-isomorphism $T_* : \K(X) \to \K(Y)$ by $T_* (\Theta_{x,y}^X) = \Theta_{T(x),T(y)}^Y$.  Also if $a \in A$ and $\phi_X(a) \in \K(X)$, then $T_*(\phi_X(a)) = \phi_Y(a)$.
\end{lemma}

\begin{proof}
It is straightforward to check that defining $T_* : \K(X) \to \K(Y)$ by $T_* (\Theta_{x,y}^X) := \Theta_{T(x),T(y)}^Y$ and extending linearly gives a $*$-homomorphism.  Since $(T^{-1})_*$ is an inverse for this $*$-homomorphism $T_*$ is a $*$-isomorphism between $C^*$-algebras.  Furthermore, if $\phi_X(a) \in \K(X)$, then $$\phi_X(a) = \lim_n \sum_{k=1}^{N_n} \Theta^X_{x_{n,k},y_{n,k}}.$$  But then for any $y \in Y$ we may let $x := T^{-1}(y)$ and we have
\begin{align*}
T_*(\phi_X(a)) y &= T_*(\phi_X(a))T(x) = T_*(  \lim_n \sum_{k=1}^{N_n} \Theta^X_{x_{n,k},y_{n,k}}) T(x) \\
&=  \lim_n \sum_{k=1}^{N_n} \Theta^Y_{T(x_{n,k}),T(y_{n,k})} T(x) = \lim_n \sum_{k=1}^{N_n} T(x_{n,k}) \langle T(y_{n,k}), T(x) \rangle_Y \\
&= \lim_n \sum_{k=1}^{N_n} T(x_{n,k}) \langle y_{n,k}, x \rangle_X = \lim_n \sum_{k=1}^{N_n} T(x_{n,k} \langle y_{n,k}, x \rangle_X) \\
&= T( \lim_n \sum_{k=1}^{N_n} x_{n,k} \langle y_{n,k}, x \rangle_X) = T(  \lim_n \sum_{k=1}^{N_n} \Theta^X_{x_{n,k},y_{n,k}} (x)) \\
&= T(\phi_X(a) x) = \phi_Y(a)T(x) = \phi_Y(a) y
\end{align*}
so that $T_*(\phi_X(a)) = \phi_Y(a)$.
\end{proof}

\subsection{Essential $C^*$-correspondences}

\begin{definition}
A $C^*$-correspondence $X$ from $A$ to $B$ is said to be \emph{essential} if $\clspan \{ \phi_X(a)x : a \in A \text{ and } x \in X \} = X$.
\end{definition}

It can be shown that $X$ is essential if and only if whenever $\{ e_\lambda \}_{\lambda \in \Lambda}$ is an approximate unit for $A$, then $\lim_\lambda \phi_X(e_\lambda)x = x$ for all $x \in X$.

\begin{definition}
If $X$ is a $C^*$-correspondence from $A$ to $B$, the \emph{essential subspace} of $X$ is defined to be $$X_\text{ess} := \clspan \{ \phi_X(a)x : a \in A \text{ and } x \in X \}.$$
\end{definition}

Notice that $X_\text{ess}$ is closed under addition and right multiplication by elements of $B$.  Thus $X_\text{ess}$ is a right Hilbert $B$-module with the inner product that it inherits from $X$.   In addition, if $a \in A$, then $\phi_X(a)|_{X_\text{ess}}$ takes values in $X_\text{ess}$ and hence $\phi_X(a)$ restricts to an element in $\Li (X_\text{ess})$.  Therefore, defining $\phi_{X_\text{ess}}(a) := \phi_X(a)|_{X_\text{ess}}$ makes $X_\text{ess}$ into a $C^*$-correspondence from $A$ to $B$.

\subsection{$C^*$-algebras associated with $C^*$-correspondences}

\begin{definition} \label{repn-def}
If $X$ is a $C^*$-correspondence over $A$, then a \emph{representation} of
$X$ into a $C^*$-algebra $B$ is a pair $(t, \pi)$ consisting of a linear map $t : X \to B$ and a
$*$-homomorphism $\pi : A \to B$ satisfying
\begin{enumerate}
\item[(i)] $t(\xi)^* t(\eta) = \pi(\langle \xi, \eta \rangle_X)$
\item[(ii)] $t(\phi_X(a)\xi) = \pi(a) t(\xi)$
\item[(iii)] $t(\xi a) = t(\xi) \pi(a)$
\end{enumerate}
for all $\xi, \eta \in X$ and $a \in A$.  We often write $(t, \pi) : (X,A) \to B$ in this situation.

Note that Condition (iii) follows from Condition (i) due to the equation $$\|
t(\xi) \pi(a) - t(\xi a) \|^2 = \| (t(\xi) \pi(a) - t(\xi a) )^* (t(\xi)
\pi(a) - t(\xi a)) \| = 0.$$  If $(t, \pi) : (X,A) \to B$ is a representation of $X$ into
a $C^*$-algebra $B$, we let $C^*(t, \pi)$ denote the $C^*$-subalgebra of
$B$ generated by $t(X) \cup \pi(A)$.
\end{definition}

\begin{definition}
A representation $(t, \pi) : (X,A) \to B$ is said to be \emph{injective} if $\pi$ is
injective.  Note that in this case $t$ will also be isometric since $$\|
t(\xi) \|^2 = \| t(\xi)^* t(\xi) \| = \| \pi (\langle \xi, \xi \rangle_X) \| =
\| \langle \xi, \xi \rangle_A \| = \| \xi \|^2.$$
\end{definition}

\begin{definition}[The Toeplitz Algebra of a $C^*$-correspondence]
Given a $C^*$-correspondence $X$ over a $C^*$-algebra $A$, there is a $C^*$-algebra $\TX$ and a representation $(\overline{t}_X, \overline{{\pi}}_X) :(X,A) \to \TX$ that is universal in the following sense: 
\begin{enumerate}
\item $\TX$ is generated as a $C^*$-algebra by $\overline{t}_X(X) \cup \overline{\pi}_X(A)$; and 
\item Given any representation $(t, \pi) :(X,A) \to B$ of $X$ into a $C^*$-algebra $B$, there exists a $*$-homomorphism of ${\rho}_{(t, \pi)} : \TX \to B$, such that $t={{\rho}_{(t, \pi)}}\circ{\overline{t}_X}$ and $\pi={{\rho}_{(t, \pi)}}\circ{\overline{\pi}}_X$.  
\end{enumerate}
The $C^*$-algebra $\TX$ and the representation $(\overline{t}_X, {\overline{\pi}}_X)$ exist (see \cite{FR}, for example) and are unique up to an obvious notion of isomorphism.  We call $\TX$ \emph{the Toeplitz algebra of the $C^*$-correspondence $X$}, and we call $(\overline{t}_X, {\overline{\pi}}_X)$ \emph{a universal representation} of $X$ in $\TX$.
\end{definition}

\begin{definition}
For a representation $(t, \pi) : (X,A) \to B$ of $X$ into a $C^*$-algebra $B$ there
exists a $*$-homomorphism $\psi_t : \K (X) \to B$ with the property that
$$\psi_t (\Theta_{\xi,\eta}) = t(\xi) t(\eta)^*.$$  
See \cite[p.~202]{Pim}, \cite[Lemma~2.2]{KPW}, and \cite[Remark~1.7]{FR} for details on the existence
of this $*$-homomorphism.  (We warn the reader that our map $\psi_t$ is denoted by $\pi^{(1)}$ in much of the literature, and by $\rho^{(t,\pi)} = \rho^{(\psi,\pi)}$ in \cite{FR}.  We have chosen to use $\psi_t$ in order to follow the conventions of Katsura in \cite{Kat3, Kat4, Kat5} and since the map depends only on $t$ and not on $\pi$.)  It is shown in \cite[Lemma~2.4]{Kat4} that if $(t,\pi)$ is an injective representation, then $\psi_t$ is injective as well.
\end{definition}

\begin{definition}
For an ideal $I$ in a $C^*$-algebra $A$ we define $$I^\perp := \{ a \in A :
ab=0 \text{ for all } b \in I \}.$$  If $X$ is a $C^*$-correspondence over
$A$, we define an ideal $J(X)$ of $A$ by $J(X) := \phi_X^{-1}(\K(X))$.  We also
define an ideal $J_X$ of $A$ by $$J_X := J(X) \cap (\ker \phi_X)^\perp.$$  Note
that $J_X = J(X)$ when $\phi_X$ is injective, and that $J_X$ is the maximal
ideal on which the restriction of $\phi$ is an injection into $\K(X)$.
\end{definition}

\begin{definition}
If $X$ is a $C^*$-correspondence over $A$, then a representation $(t, \pi) : (X,A) \to B$ of $X$ into a $C^*$-algebra $B$ is said to be \emph{coisometric} if $$\psi_t (\phi_X(a)) = \pi(a) \qquad \text{for all $a \in J_X$.}$$
\end{definition}

\begin{definition}[The Cuntz-Pimsner Algebra of a $C^*$-correspondence]
Given a $C^*$-correspondence $X$ over a $C^*$-algebra $A$, there is a $C^*$-algebra $\OX$ and a coisometric representation $(t_X, {\pi}_X) :(X,A) \to \OX$ that is universal in the following sense: 
\begin{enumerate}
\item $\OX$ is generated as a $C^*$-algebra by $t_X(X) \cup {\pi}_X(A)$; and 
\item Given any coisometric representation $(t, \pi) :(X,A) \to B$ of $X$ into a $C^*$-algebra $B$, there exists a $*$-homomorphism of ${\rho}_{(t, \pi)} : \OX \to B$, such that $t={{\rho}_{(t, \pi)}}\circ{t_X}$ and $\pi={{\rho}_{(t, \pi)}}\circ{\pi}_X$.  
\end{enumerate}
The $C^*$-algebra $\OX$ and the representation $(t_X, \pi_X)$ exist (see \cite[\S4]{Kat4}) and are unique up to an obvious notion of isomorphism.  We call $\OX$ \emph{the Cuntz-Pimsner algebra of the $C^*$-correspondence $X$}, and we call $(t_X, {\pi}_X)$ \emph{a universal coisometric representation} of $X$ in $\OX$.  We also mention that any universal coisometric representation $(t_X, {\pi}_X)$ is injective.
\end{definition}

\subsection{The Gauge Action on Cuntz-Pimsner Algebras}

\begin{definition}
If $X$ is a $C^*$-correspondence over $A$, we say that a representation $(t, \pi) : (X,A) \to B$ of $X$ into a $C^*$-algebra $B$ \emph{admits a gauge action} if for each $z \in \T$ there is a $*$-homomorphism $\beta_z : C^*(t,\pi) \to C^*(t,\pi)$ such that $\beta_z(t(\xi)) = zt(\xi)$ for all $\xi \in X$ and $\beta_z (\pi(a)) = \pi(a)$ for all $a \in A$.

It is a consequence of this definition that $\beta_z$ is actually an automorphism with $\beta_z^{-1} = \beta_{\overline{z}}$, and that the map $\beta : \T \to \aut C^*(t,\pi)$ given by $z \mapsto \beta_z$ is strongly continuous.
\end{definition}

\begin{definition}
If $X$ is a $C^*$-correspondence over $A$, and $(t_X, \pi_X) : (X,A) \to \OX$ is the universal coisometric representation of $X$ into $\OX$, then the universal property of $\OX$ implies that $(t_X, \pi_X)$ admits a gauge action, which we denote by $\gamma : \T \to \OX$.  We refer to $\gamma$ as \emph{the canonical gauge action on $\OX$}.
\end{definition}

\begin{theorem}[The Gauge-Invariant Uniqueness Theorem for Cuntz-Pimsner Algebras]
Let $X$ be a $C^*$-correspondence over a $C^*$-algebra $A$, and let $(t, \pi) : (X,A) \to B$ be a coisometric representation of $X$ into a $C^*$-algebra $B$.  If $C^*(t,\pi)$ is the $C^*$-subalgebra of $B$ generated by the images, $t(X)$ and $\pi(A)$, then the induced $*$-homomorphism $\rho_{(t,\pi)} : \OX \to C^*(t,\pi)$ is an isomorphism if and only if $(t,\pi)$ is injective and admits a gauge action.
\end{theorem}
 
The Gauge-Invariant Uniqueness Theorem is proven in \cite[\S6]{Kat4}.  It is one of our most important tools for constructing isomorphisms of Cuntz-Pimsner algebras.

\subsection{Tensor products of $C^*$-correspondences}

If ${}_AX_B$ is a $C^*$-correspondence from $A$ to $B$, and if ${}_BY_C$ is a $C^*$-correspondence from $B$ to $C$, then we may form a correspondence $X \otimes_B Y$ from $A$ to $C$, called the \emph{internal tensor product} (sometimes also called the \emph{interior tensor product}) as follows:   We first regard $Y$ as a left $B$-module and form the algebraic tensor product $X \odot Y$.  We then let $N$ be the subspace generated by $$\{ xb \odot y - x \odot \phi_Y(b) y : x \in X, y \in Y, \text{ and } b \in B \}$$ and form the balanced tensor product $X \odot_B Y := (X \odot Y) / N$.  If $x \odot y$ is an elementary tensor in $X \odot Y$, we let $x \odot_B y$ denote its equivalence class in $X \odot_B Y$.  We give $X \odot Y$ a right $C$ action by defining $(x \odot_B y) c := x \odot_B yc$, a left $A$ action by defining $\phi_{X \odot_B Y} (a) (x \odot_B y) = \phi_X(a)x \odot_B y$, and a $C$-valued inner product by defining $$\langle x_1 \odot_B y_1, x_2 \odot_B y_2 \rangle_{X \odot_B Y} := \langle y_1, \phi_Y(\langle x_1, x_2 \rangle_X) y_2 \rangle_Y.$$  These formulae are well-defined and continuous on all of $X \odot_B Y$ and make $X \odot_B Y$ into a pre-$C^*$-correspondence from $A$ to $C$.    (In particular, the subspace $\{ z \in X \odot_B Y : \langle z, z \rangle_{X \odot_B Y} = 0 \}$ is equal to $N$  \cite[p.40]{Lan} so that the inner product on $X \odot_B Y$ is nondegenerate.)  We then define $X \otimes_B Y$ to be the completion of $X \odot_B Y$ with respect to the norm coming from the above inner product, and we let $x \otimes y$ denote the equivalence class of $x \odot y \in X \odot Y$.  We mention that $X \otimes_B Y = \clspan \{ x \otimes y : x \in X \text{ and } y \in Y \}$.  If $T \in \Li(X)$ and $S \in \Li(Y)$ with $\phi_Y(b) S(y) = S (\phi_Y(b)y)$ for all $b \in B$ and $y \in Y$, then one can show that there exists an operator $T \otimes_A S \in \Li (X \otimes_A Y)$ with $(T \otimes_A S) (x \otimes y) = T(x) \otimes S(y)$.

Note that if $T \in \K(X \otimes_B Y)$, then the linearity of the inner product, and the fact that $\operatorname{span} \{ x \otimes y : x \in X \text{ and } y \in Y \}$ is dense in $X \otimes_B Y$ allows us to write $T$ as the limit of finite sums of elements $\Theta^{X \otimes_B Y}_{x \otimes y, z \otimes w}$; that is, the subscripts in the generalized rank one operators may be chosen to be elementary tensors.

If $X$ is a $C^*$-correspondence over $A$, then we may form the tensor product of $X$ with itself.  For $n \geq 1$ we let $X^{\otimes n}$ denote the tensor product $X \otimes_A \ldots \otimes_A X$ of $n$ copies of $X$.  We then have that $$X^{\otimes n} = \clspan \{\xi_1 \otimes \ldots \otimes \xi_n : \xi_1, \ldots, \xi_n \in X \}.$$  If $(t,\pi) : (X,A) \to B$ is a representation of $X$ into a $C^*$-algebra $B$, then it is straightforward to show that there exists a linear map $t^n : X^{\otimes n} \to B$ defined by $t^n(\xi_1 \otimes \ldots \otimes \xi_n) = t(\xi_1) \ldots t(\xi_n)$, and that $(t^n,\pi) : (X^{\otimes n}, A) \to B$ is a representation of $X^{\otimes n}$ into $B$.  In particular, there exists a $*$-homomorphism $\psi_{t^n} : \K(X^{\otimes n}) \to B$ with $$\psi_{t^n} ( \Theta_{\xi, \eta}^{X^{\otimes n}} ) = t^n(\xi)t^n(\eta)^* \qquad \text{ for } \xi, \eta \in X^{\otimes n}.$$  For $n=0$ we define $X^{\otimes 0} =A$, and we take $t^0 := \pi$.  It can also be shown (see \cite[Proposition~2.7]{Kat4}) that $$C^*(t,\pi) = \clspan \{ t^n(\xi)t^m(\eta)^* : \xi \in X^{\otimes n}, \eta \in X^{\otimes m}, \text{ and } n,m \in \N \}.$$

The following proposition was proven in \cite[Lemma~4.5]{FMR} and will be useful for us in our analysis.

\begin{proposition} \label{FMR-prop}
Let $X$ be a $C^*$-correspondence over $A$ and suppose that the left action $\phi_X : A \to \Li(X)$ is injective.  Also let $n \geq 1$.  If $S \in \Li(X^{\otimes n})$ and $S \otimes \textrm{Id} \in \K (X^{\otimes (n+1)})$, then $S \in \K (X^{\otimes n})$ and $\psi_{t^{n+1}} (S \otimes_A \textrm{Id}) = \psi_{t^n} (S)$.
\end{proposition}

\section{Elementary Strong Shift Equivalence of Regular $C^*$-correspondences} \label{ESSERCC}

\begin{definition}
We say that a $C^*$-correspondence ${}_AX_B$ from $A$ to $B$ is \emph{regular} if the left action $\phi_X : A \to \Li(X)$ is injective and $\im \phi_X \subseteq \K(X)$.  Note that if $X$ is regular, then $J_X = J(X) = A$.
\end{definition}

\begin{definition}
Let ${}_AE_A$ be a $C^*$-correspondence over $A$ and let ${}_BF_B$ be a $C^*$-correspondence over $B$.  We say that $E$ is \emph{elementary strong shift equivalent} to $F$ if there exists a $C^*$-correspondence ${}_AR_B$ from $A$ to $B$ and a $C^*$-correspondence ${}_BS_A$ from $B$ to $A$ such that $$E \cong R \otimes_B S \qquad \text{ and } \qquad F \cong S \otimes_A R.$$
\end{definition}

We shall spend the remainder of this section proving that if $E$ and $F$ are essential, regular $C^*$-correspondences that are elementary strong shift equivalent, then $\mathcal{O}_E$ is Morita equivalent to $\mathcal{O}_F$.

\begin{definition} \label{def-bipart-infl}
Let ${}_AR_B$ be a $C^*$-correspondence from $A$ to $B$ and let ${}_BS_A$ be a $C^*$-correspondence  from $B$ to $A$.  The \emph{bipartite inflation of $S$ by $R$} is a $C^*$-correspondence $X$ over $A \oplus B$ defined in the following way:  We let $$X = S \oplus R$$ be a right Hilbert $A \oplus B$-module as in Remark~\ref{direct-sum-right}, and we make $X$ into a $C^*$-correspondence $X$ over $A \oplus B$ by defining $$\phi_X(a,b) (s,r) := (\phi_S(b)s, \phi_R(a)r).$$
Note that the order of $S$ and $R$ are relevant, and in particular, the bipartite inflation of $S$ by $R$ is not equal to the bipartite inflation of $R$ by $S$.
\end{definition}

Throughout this section fix $C^*$-algebras $A$ and $B$, and fix a $C^*$-correspondence ${}_AR_B$ from $A$ to $B$ and a $C^*$-correspondence ${}_BS_A$ from $B$ to $A$.  We shall set
\begin{equation} \label{E-F-RS}
E = R \otimes_B S \qquad \text{ and } \qquad F = S \otimes_A R
\end{equation}
so that $E$, which is a $C^*$-correspondence over $A$, is elementary strong shift equivalent to $F$, which is a $C^*$-correspondence over $B$.  We shall also let $X = S \oplus R$ denote the bipartite inflation of $S$ by $R$, which is a $C^*$-correspondence over $A \oplus B$.

\begin{proposition} \label{iso-X2}
With the notation above, we have $X^{\otimes 2} \cong E \oplus F$ as $C^*$-correspondences.  Furthermore, there exists an isomorphism $T : X^{\otimes 2} \to E \oplus F$ with the property that $T((s,r) \otimes (s', r')) = ( r \otimes s', s \otimes r')$.
\end{proposition}

\begin{proof}
We begin by defining a balanced bilinear map $T_0 : X \oplus X \to E \oplus F$ by $T_0((s,r),(s',r')) = (r \otimes s', s \otimes r')$.  Then $T_0$ induces a linear map $T : X \odot_{A \oplus B} X \to E \oplus F$, and for any pair of elementary tensors $(s_1,r_1) \odot_{A \oplus B} (s_1', r_1'), (s_2,r_2) \odot_{A \oplus B} (s_2', r_2') \in X \odot_{A \oplus B} X$ we have
\begin{align*}
&\langle T ((s_1,r_1) \odot_{A \oplus B} (s_1', r_1')) , T((s_2,r_2) \odot_{A \oplus B} (s_2', r_2')) \rangle_{E \oplus F} \\
=& \langle (r_1 \otimes s_1', s_1 \otimes r_1'), (r_2 \otimes s_2', s_2 \otimes r_2') \rangle_{E \oplus F} \\
=& (\langle r_1 \otimes s_1', r_2 \otimes s_2' \rangle_E, \langle s_1 \otimes r_1', s_2 \otimes r_2'\rangle_F) \\
=& (\langle s_1', \phi_R(\langle r_1,r_2 \rangle_R) s_2' \rangle_S, \langle r_1', \phi_S(\langle s_1, s_2 \rangle_S) r_2'\rangle_R) \\
=& \langle (s_1', r_1') , ( \phi_R(\langle r_1,r_2 \rangle_R) s_2', \phi_S(\langle s_1, s_2 \rangle_S) r_2' ) \rangle_X \\
=& \langle (s_1', r_1') , \phi_X( \langle s_1, s_2 \rangle_S, \langle r_1,r_2 \rangle_R) (s_2',  r_2' ) \rangle_X \\
=& \langle (s_1,r_1) \odot_{A \oplus B} (s_1', r_1') , (s_2,r_2) \odot_{A \oplus B} (s_2', r_2') \rangle_{X^{\otimes 2}}.
\end{align*}
Because the inner product is bilinear, this shows that $\langle T(x), T(y) \rangle_{E \oplus F} = \langle x, y \rangle_{X^{\otimes 2}}$ for all $x, y \in X \odot_{A \oplus B} X$.  Thus $T$ is bounded and extends to a map $T : X^{\otimes 2} \to E \oplus F$, which preserves inner products.  Furthermore, it is straightforward to check that $T(\phi_{X^{\otimes 2}}(a,b) z) = \phi_{E \oplus F}(a,b) T(z)$ for all $(a,b) \in A \oplus B$ and $z \in X^{\otimes 2}$.  Finally, since $T$ preserves inner products we have that $T$ is injective, and since $\clspan \{ r \otimes s' : r \in R \text{ and } s' \in S \} = R \otimes_B S = E$ and $\clspan \{ s \otimes r' : s \in S \text{ and } r' \in R \} = S \otimes_A R = F$ we see that $T$ is surjective.   Thus $T$ is an isomorphism of $C^*$-correspondences.
\end{proof}

\begin{lemma} \label{E-F-inj-X-inj}
If $\phi_{R \otimes_B S}$ and $\phi_{S \otimes_A R}$ are injective, then $\phi_X$ is injective.
\end{lemma}

\begin{proof}
Since $\phi_E(a) = \phi_R(a) \otimes \textrm{Id}$ for all $a \in A$, we see that if $\phi_R(a) = 0$, then $\phi_{R \otimes_B S}(a) = 0$ and the injectivity of $\phi_{R \otimes_B S}$ implies that $a = 0$.  Thus $\phi_R$ is injective.  Similarly, the injectivity of $\phi_{S \otimes_A R}$ implies that $\phi_S$ is injective.  Because $X = S \oplus R$ and $\phi_X(a,b) = \phi_S(b) \oplus \phi_R(a)$ for all $(a,b) \in A \oplus B$, we have that $\phi_X$ is injective.
\end{proof}

\begin{lemma} \label{phi-RS-X2}
If $a \in J_{R \otimes_B S}$ and $\phi_{R \otimes_B S}(a) = \lim_n \sum_{k=1}^{N_n} \Theta^{R \otimes_B S}_{r_{n,k} \otimes s_{n,k}, r_{n,k}' \otimes s_{n,k}'}$, then $\phi_{X^{\otimes 2}}(a,0) = \lim_n \sum_{k=1}^{N_n} \Theta^{X^{\otimes 2}}_{(0,r_{n,k}) \otimes (s_{n,k},0), (0,r_{n,k}') \otimes (s_{n,k}',0)}$.
\end{lemma}

\begin{proof}
Since $E = R \otimes_B S$ and $\phi_{E \oplus F} = \phi_E \oplus \phi_F$, we see that $\phi_{E \oplus F} (a,0) = \lim_n \sum_{k=1}^{N_n} \Theta^{E \oplus F}_{(r_{n,k} \otimes s_{n,k},0), (r_{n,k}' \otimes s_{n,k}',0)}$.  Using the isomorphism, $((r \otimes s'),(s \otimes r') \mapsto (s,r) \otimes (s',r')$ from $E \oplus F$ to $X^{\otimes 2}$ established in Proposition~\ref{iso-X2}, we may then apply Lemma~\ref{iso-phi-translate} to obtain the result.
\end{proof}

\begin{lemma} \label{t-pi-X-def}
Suppose $\phi_{R \otimes_B S}$ and $\phi_{S \otimes_A R}$ are injective.  Let $(t_X, \pi_X) : (X,A \oplus B) \to \OX$ be a universal coisometric representation of $X$ into $\OX$.  Then there exists a coisometric representation $(t,\pi) :(R \otimes_B S , A) \to \OX$ with $$t(r \otimes s) = t_X(0,r) t_X(s,0) \qquad \text{ and } \qquad \pi(a) = \pi_X(a,0)$$ for all $r \otimes s \in R \otimes_B S$ and for all $a \in A$.
\end{lemma}

\begin{proof}
Begin by defining $\pi : A \to \OX$ by $\pi(a) := \pi_X(a,0)$.  Since $(r,s) \mapsto (0,r) \otimes (s,0)$ is a balanced bilinear map from $R \oplus S$ to $X^{\otimes 2}$, it induces a linear map $T : R \odot_B S \to X^{\otimes 2}$.  If $r \odot_B s$ and $r' \odot_B s'$ are elementary tensors in $R \odot_B S$, then 
\begin{align*}
\langle T(r \odot_B s), T(r' \odot_B s') \rangle_{X^{\otimes 2}} &= \langle (0,r) \otimes (s,0) , (0,r') \otimes (s',0) \rangle_{X^{\otimes 2}} \\
&=  \langle(s,0) , \phi_X(\langle (0,r), (0,r') \rangle_X) (s',0) \rangle_X \\
&= (\langle s, \phi_R(\langle r, r' \rangle_R) s', 0 ) \\
&= (\langle r \otimes s, r' \otimes s' \rangle_{R \otimes S}, 0).
\end{align*}
By the bilinearity of the inner product, it follows that $\langle T(z), T(w) \rangle_{X^{\otimes 2}} = (\langle z, w \rangle_{R \otimes_B S}, 0)$ for all $z,w \in R \odot_B S$.  Thus $|T(z)| = | \langle T(z), T(z) \rangle_{X^{\otimes 2}}|^{1/2} = |(\langle z, z \rangle_{R \otimes_B S}, 0)|^{1/2} = |z|$ so $T$ is isometric and extends to a map $T : R \otimes_B S \to X^{\otimes 2}$.  If we let $t := t_X^2 \circ T$, then $t:R \otimes_B S \to \OX$ and $t(r \otimes s) := t_X(0,r) t_X(s,0).$

We shall now show that $(t, \pi)$ is a representation of $R \otimes_B S$ into $\OX$.  If $a \in A$ and $r \otimes s$ and $r' \otimes s'$ are elementary tensors in $R \otimes_B S$, then 
\begin{align*}
t(r \otimes s)^* t(r' \otimes s') &= t_X(s,0)^* t_X(0,r)^* t_X(0,r') t_X(s',0) \\
&= t_X(s,0)^* \pi_X(\langle (0,r), (0,r') \rangle_X) t_X(s',0) \\
&=  t_X(s,0)^* \pi_X((0,\langle r, r' \rangle_R)) t_X(s',0) \\
&= \pi_X (\langle s, \phi_S(\langle r, r' \rangle_R)s'\rangle_S, 0) \\
&= \pi(\langle r \otimes s, r' \otimes s' \rangle_{R \otimes_B S}).
\end{align*}
and
\begin{align*}
t(\phi_{R \otimes_B S}(a) (r \otimes s)) &= t_X(0,\phi_R(a)r) t_X(s,0) = t_X (\phi_X(a,0)(0,r)) t_X(s,0) \\
&=  \pi_X(a,0) t_X(0,r) t_X(s,0) = \pi(a)t(r \otimes s).
\end{align*}
Because of linearity and the fact that elementary tensors span a dense subset of $X^{\otimes 2}$, the above two equations show that Condition~(i) and Condition~(ii) of Definition~\ref{repn-def} hold and $(t, \pi)$ is a representation of $R \otimes_B S$.

To see that $(t, \pi)$ is coisometric, let $a \in J_{R \otimes_B S}$ and write $\phi_{R \otimes_B S}(a) = \lim_n \sum_{k=1}^{N_n} \Theta^{R \otimes_B S}_{r_{n,k} \otimes s_{n,k}, r_{n,k}' \otimes s_{n,k}'}$.  Then 
\begin{equation} \label{phi-a-compact}
\phi_{X^{\otimes 2}}(a,0) = \lim_n \sum_{k=1}^{N_n} \Theta^{X^{\otimes 2}}_{(0,r_{n,k}) \otimes (s_{n,k},0), (0,r_{n,k}') \otimes (s_{n,k}',0)}
\end{equation}
by Lemma~\ref{phi-RS-X2}.  Hence
\begin{align*}
&\psi_t(\phi_{R \otimes_B S}(a)) \\
= &\lim_n \sum_{k=1}^{N_n} t(r_{n,k} \otimes s_{n,k}) t(r'_{n,k} \otimes s'_{n,k})^* \\
= &\lim_n \sum_{k=1}^{N_n} t_X(0, r_{n,k}) t_X (s_{n,k},0) (t_X(0, r'_{n,k}) t_X (s'_{n,k},0))^* \\
= &\lim_n \sum_{k=1}^{N_n} t_X^{2}((0, r_{n,k}) \otimes (s_{n,k},0)) t_X^{2}((0, r'_{n,k}) \otimes (s'_{n,k},0))^* \\
= & \psi_{t_{X^{2}}} \left( \lim_n \sum_{k=1}^{N_n} \Theta^{X^{\otimes 2}}_{(0,r_{n,k}) \otimes (s_{n,k},0), (0,r_{n,k}') \otimes (s_{n,k}',0)} \right) \\
= & \psi_{t_{X^{2}}} (\phi_{X^{\otimes 2}}(a,0)) \\
= & \psi_{t_{X^{2}}} (\phi_X(a,0) \otimes \textrm{Id}).
\end{align*}
From Lemma~\ref{E-F-inj-X-inj} it follows that $\phi_X$ is injective.  Also, $\phi_X(a,0) \otimes \textrm{Id} = \phi_{X^{\otimes 2}}(a,0) \in \K(X)$ from (\ref{phi-a-compact}).  Thus the above equation and Proposition~\ref{FMR-prop} show that $\psi_t(\phi_{R \otimes_B S}(a)) = \psi_{t_{X^{2}}} (\phi_X(a,0) \otimes \textrm{Id}) = \psi_{t_X}(\phi_X(a,0)) = \pi_X(a,0) = \pi(a)$ and $(t,\pi)$ is coisometric.
\end{proof}

\begin{lemma} \label{gauge-act-corner}
Suppose $\phi_{R \otimes_B S}$ and $\phi_{S \otimes_A R}$ are injective.  If $(t,\pi) :(R \otimes_B S, A) \to \OX$ is the coisometric representation defined in Lemma~\ref{t-pi-X-def}, then $(t,\pi)$ admits a gauge action.
\end{lemma}

\begin{proof}
Let $(t_X,\pi_X) : (X,A) \to \OX$ be a universal coisometric representation of $X$ into $\OX$.  For $z \in \T$ define $(t_z,\pi_z) :(X,A) \to \OX$ by $\pi_z = \pi_X$ and $t_z(s,r) := t_X(sz, r)$ for $(s,r) \in X = S \oplus R$.  Then $t_z$ is a linear map, and we see that $(t_z, \pi_z)$ is a representation since:
\begin{align*}
t_z(s,r)^*t_z(s',r') &= t_X(sz,r)^*t_X(s'z,r) = \pi_X(\langle (sz, r), (s'z,r') \rangle_X) \\
&= \pi_X(\langle sz,s'z \rangle_S, \langle r, r' \rangle_R) = \pi_X (\langle s,s' \rangle_S, \langle r, r' \rangle_R) \\
&= \pi_X(\langle (s,r), (s',r') \rangle_X) = \pi_z (\langle (s,r), (s',r') \rangle_X)
\end{align*}
and
\begin{align*}
t_z(\phi_X(a,b)(s,r)) &= t_z(\phi_S(b)s, \phi_R(a)r) = t_X(\phi_S(b) sz, \phi_R(a) r) \\
&= t_X(\phi_X(a,b) (sz, r)) = \pi_X(a,b)t_X(sz,r) = \pi_z(a,b)t_z(s,r).
\end{align*}
We shall also show that $(t,\pi)$ is coisometric.  If $(a,b) \in J_X$, then $$\phi_X(a,b) = \lim_n \sum_{k=1}^{N_n} \Theta^X_{(s_{n,k}.r_{n,k}),(s'_{n,k},r'_{n,k})}.$$  Since $X = S \oplus R$ and $\phi_X(a,b) = \phi_S(b) \oplus \phi_R(a)$ it follows that $\phi_S(b) = \lim_n \sum_{k=1}^{N_n} \Theta^S_{s_{n,k},s'_{n,k}}$ and $\phi_R(a) = \lim_n \sum_{k=1}^{N_n} \Theta^R_{r_{n,k},r'_{n,k}}$.  Thus
\begin{align*}
& \ \psi_{t_z} (\phi_X(a,b)) \\
=& \ \psi_{t_z}(\phi_X(a,0)) + \psi_{t_z}(\phi_X(0,b)) \\
=&  \ \psi_{t_z}(\lim_n \sum_{k=1}^{N_n} \Theta^X_{(0,r_{n,k}),(0,r'_{n,k})}) + \psi_{t_z}(\lim_n \sum_{k=1}^{N_n} \Theta^X_{(s_{n,k},0),(s'_{n,k},0)}) \\
=& \ \lim_n \sum_{k=1}^{N_n} t_z(0,r_{n,k}) t_z(0,r'_{n,k})^* + \lim_n \sum_{k=1}^{N_n} t_z(s_{n,k},0) t_z(s'_{n,k},0)^* \\
=& \ \lim_n \sum_{k=1}^{N_n} t_X(0,r_{n,k}) t_X(0,r'_{n,k})^* + \lim_n \sum_{k=1}^{N_n} t_X(s_{n,k},0)z (t_X(s'_{n,k},0)z)^* \\
=& \ \lim_n \sum_{k=1}^{N_n} t_X(0,r_{n,k}) t_X(0,r'_{n,k})^* + \lim_n \sum_{k=1}^{N_n} t_X(s_{n,k},0) t_X(s'_{n,k},0)^* \\
=& \ \psi_{t_X} \left( \lim_n \sum_{k=1}^{N_n} \Theta^X_{(0,r_{n,k}),(0,r'_{n,k})} + \lim_n \sum_{k=1}^{N_n} \Theta^X_{(s_{n,k},0),(s'_{n,k},0)} \right) \\
=& \ \psi_{t_X} (\phi_X(a,0) + \phi_X(0,b)) \\
=& \ \psi_{t_X} (\phi_X(a,b)) \\
=& \ \pi_X(a,b) \\
=& \ \pi_z(a,b).
\end{align*}
Since $(t_z,\pi_z)$ is a coisometric representation, it induces a $*$-homomorphism $\beta_z := \rho_{(t_z,\pi_z)} : \OX \to \OX$ with $\beta_z \circ t_X = t_z$ and $\beta_z \circ \pi_X = \pi_z$.  

If $(t, \pi) : (R \otimes_B S, A) \to \OX$ is the coisometric representation defined in Lemma~\ref{t-pi-X-def}, then for any elementary tensor $r \otimes s$ in $R \otimes_B S$ we have $\beta_z(t(r \otimes s)) = \beta_z(t_X(0,r)t_X(s,0)) = t_z(0,r)t_z(s,0) = t_X(0,r)t_X(sz,0) = zt_X(0,r)t_X(s,0) = z t(r \otimes s)$.  Since the elementary tensors span a dense subset of $R \otimes_B S$ it follows that $\beta_z(t(\xi)) = t(\xi)$ for all $\xi \in R \otimes_B S$.  In addition, $\beta_z(\pi(a)) = \beta_z(\pi_X(a,0)) = \pi_z(a,0) = \pi_X(a,0) = \pi(a)$ for all $a \in A$.  Thus $(t,\pi)$ admits a gauge action.
\end{proof}

\begin{lemma} \label{t-pi-prime-X-def}
Suppose $\phi_{R \otimes_B S}$ and $\phi_{S \otimes_A R}$ are injective.  Let $(t_X, \pi_X) : (X,A \oplus B) \to \OX$ be a universal coisometric representation of $X$ into $\OX$.  Then there exists a coisometric representation $(t',\pi') :(S \otimes_A R , B) \to \OX$ with $$t'(s \otimes r) = t_X(s,0) t_X(0,r) \qquad \text{ and } \qquad \pi'(b) = \pi_X(0,b)$$ for all $s \otimes r \in S \otimes_A R$ and for all $b \in B$.  Furthermore, $(t',\pi')$ admits a gauge action.
\end{lemma}

\begin{proof}
Let $X = S \oplus R$ be the bipartite inflation of $S$ by $R$, and let $(t_X, \pi_X) : (X, A \oplus B) \to \OX$ be a universal coisometric representation of $X$ into $\OX$.  Also let  $Y = R \oplus S$ be the bipartite inflation of $R$ by $S$ so that, in particular, $Y$ is a right Hilbert $B \oplus A$-module as in Remark~\ref{direct-sum-right}, and $Y$ made into a $C^*$ correspondence over $B \oplus A$ by defining $$\phi_Y(b,a) (r,s) := (\phi_R(a)r, \phi_S(b)s).$$  We may define a representation $(t_Y,\pi_Y) :(Y, B \oplus A) \to \OX$ by letting $t_Y( r,s) := t_X(s,r)$ for $(r,s) \in Y$, and $\pi_Y(b,a) := \pi_X(a,b)$ for $(b,a) \in B \oplus A$.

It is straightforward to verify that $(t_Y, \pi_Y)$ is a coisometric representation of $Y$ into $\OX$.  Furthermore, it is also straightforward to show that $(t_Y, \pi_Y)$ is universal.  (In particular, this implies that $\mathcal{O}_Y \cong \OX$.)

We may now apply Lemma~\ref{t-pi-X-def} (after interchanging the roles of $R$ and $S$ and the roles of $A$ and $B$ in the statement of the lemma) to conclude that there exists a coisometric representation $(t',\pi') : (S \otimes_A R, B) \to \OX$ with $t'(s \otimes r) = t_Y(0,s)t_Y(r,0)$ and $\pi'(b) = \pi_Y(b,0)$.  But then $t'(s \otimes r) = t_X(s,0) t_X(0,r)$ and $\pi'(b) = \pi_X(0,b)$.  Furthermore, Lemma~\ref{gauge-act-corner} shows that $(t',\pi')$ admits a gauge action.
\end{proof}

\begin{lemma} \label{ess-across-otimes}
Let $E$ be a $C^*$-correspondence over $A$.  Suppose $R$ is a $C^*$-correspondence from $A$ to $B$, and $S$ is a $C^*$-correspondence from $B$ to $A$ with the property that $E \cong R \otimes_B S$.  Then $$E_\text{ess} \cong R_\text{ess} \otimes_B S_\text{ess}.$$
\end{lemma}

\begin{proof}
Since $(r,s) \mapsto r \otimes s$ is a balanced bilinear mapping from $R_\text{ess} \oplus S_\text{ess}$ to $R \otimes_B S$, it induces a linear map $T : R_\text{ess} \odot_B S_\text{ess} \to R \otimes_B S$ with $T(r \odot_B s) = r \otimes s$.  Furthermore, because 
\begin{align*}
\langle T(r \odot_B s), &T( r' \odot_B s') \rangle_{R \otimes_B S} = \langle r \otimes s, r' \otimes s' \rangle_{R \otimes_B S} =  \langle s, \phi_R(\langle r,r' \rangle_R) s' \rangle_S \\
&=  \langle s, \phi_{R_\text{ess}}(\langle r,r' \rangle_{R_\text{ess}}) s' \rangle_{S_\text{ess}} =  \langle r \odot_B s, r' \odot_B s' \rangle_{R_\text{ess} \otimes_B S_\text{ess}}
\end{align*}
we see that $T$ extends to an isometric map $T: R_\text{ess} \otimes_B S_\text{ess} \to R \otimes S$.  Furthermore, if $r \in R_\text{ess}$ and $s \in S_\text{ess}$, then by the Hewitt-Cohen factorization Theorem \cite[Proposition~2.33]{RW} we may write $r = \phi_R(a)r'$.  Thus $$T(r \otimes s) = \phi_R(a)r' \otimes s = \phi_{R \otimes_B S}(a) (r' \otimes s) \in (R \otimes_B S)_\text{ess}$$ and hence $\im T \subseteq (R \otimes S)_\text{ess}$.

For the reverse inclusion, choose any such elementary tensor $\phi_R(a)r \otimes s$, and let $\{ e_\lambda \}_{\lambda \in \Lambda}$ be an approximate unit for $B$.  Then $\lim_\lambda r e_\lambda = r$ and 
\begin{align*}
\phi_R(a)r \otimes s &= \lim_\lambda \phi_R(a)re_\lambda \otimes s = \lim_\lambda \phi_R(a)r \otimes \phi_S(e_\lambda)s \\
&=  \lim_\lambda T (\phi_R(a)r \otimes \phi_S(e_\lambda)s) = T(  \lim_\lambda \phi_R(a)r \otimes \phi_S(e_\lambda)s)
\end{align*}
so that $\phi_R(a) \otimes s$ is in the image of $T$, and because the span of elementary tensors of the form $\phi_R(a)r \otimes s$ is dense in $(R \otimes_B S)_\text{ess}$ we have that $\im T = (R \otimes_B S)_\text{ess}$.  Thus $T$ is an isomorphism from $R_\text{ess} \otimes_B S_\text{ess}$ onto $(R \otimes S)_\text{ess}$, and since $E \cong R \otimes_B S$ we have $E_\text{ess} \cong R_\text{ess} \otimes_B S_\text{ess}$ .
\end{proof}

\begin{corollary}  \label{choose-RS-essential}
Let $E$ be a $C^*$-correspondence over $A$, let $F$ be a $C^*$-correspondence over $B$, and suppose that $E$ and $F$ are elementary strong shift equivalent.  If $E$ and $F$ are essential, then there exists an essential $C^*$-correspondence $R'$ from $A$ to $B$ and an essential $C^*$-correspondence $S'$ from $B$ to $A$ for which $$E \cong R' \otimes_B S' \qquad \text{ and } \qquad F \cong S' \otimes_A R'.$$
\end{corollary}

\begin{proof}
Since $E$ and $F$ are elementary strong shift equivalent there exist $R$ and $S$ with $E \cong R \otimes_B S$ and $F \cong S \otimes_A R$.  Because $E = E_\text{ess} \cong R_\text{ess} \otimes_B S_\text{ess}$ and $F = F_\text{ess} \cong S_\text{ess} \otimes_A R_\text{ess}$, we may take $R' = R_\text{ess}$ and $S' = S_\text{ess}$.
\end{proof}

\begin{lemma} \label{proj-mult-alg}
Suppose $R$ and $S$ are essential, and let $X$ be the bipartite inflation of $S$ by $R$.  If $(t_X, \pi_X) :(X,A \oplus B) \to \OX$ is a universal coisometric representation of $X$ into $\OX$, then there exist projections $P_E$ and $P_F$ in the multiplier algebra $\M (\OX)$ such that 
\begin{enumerate}
\item $P_E \ t_X(s,r) = t_X(s,0)$, 
\item $t_X(s,r) \ P_E = t_X(0,r)$, 
\item $P_E \ \pi_X(a,b) = \pi_X(a,0)$, 
\item $P_F \ t_X(s,r) = t_X(0, r)$, 
\item $t_X(s,r) \ P_F= t_X(s, 0)$, and
\item $P_F \ \pi_X(a,b) = \pi_X(0,b)$
\end{enumerate}
for all $(s,r) \in X$ and $(a,b) \in A \oplus B$.
\end{lemma}

\begin{proof}
Let $\{ e_\lambda \}_{\lambda \in \Lambda}$ be an approximate unit for $B$.  Since $S$ is essential, $\lim_\lambda \phi_S(e_\lambda)s=s$ for all $a \in A$.   For any element 
\begin{equation} \label{proj-approx}
t_X(s_1,r_1) \ldots t_X(s_n,r_n) t_X(s_m'r_m')^* \ldots t_X(s_1',r_1')^*
\end{equation}
we have that 
\begin{align*}
& \ \lim_\lambda \pi_X(0,e_\lambda) t_X(s_1,r_1) \ldots t_X(s_n,r_n) t_X(s_m'r_m')^* \ldots t_X(s_1',r_1')^* \\
= & \ \lim_\lambda t_X(\phi_S(e_\lambda) s_1, 0) t_X(s_2,r_2) \ldots t_X(s_n,r_n) t_X(s_m'r_m')^* \ldots t_X(s_1',r_1')^* \\
= &  \ t_X( s_1, 0) t_X(s_2,r_2) \ldots t_X(s_n,r_n) t_X(s_m'r_m')^* \ldots t_X(s_1',r_1')^*
\end{align*}
so this limit exists.  Because any element of $\OX$ can be approximated by a finite sum of elements of the form shown in (\ref{proj-approx}), we see that $\lim_\lambda \pi_X(0,e_\lambda) x$ exists for all $x \in \OX$.  Let us view $\OX$ as a $C^*$-correspondence over itself (see \cite[Example 2.10]{RW}).  If we define $P_E : \OX \to \OX$ by $P_E(x) := \lim_\lambda \pi_X(0,e_\lambda)x$ then we see that for any $x, y \in \OX$ we have $$y^* P_E(x) = \lim_\lambda y^*\pi_X(0,e_\lambda)x = \lim_\lambda (\pi_X(0,e_\lambda) y)^* x = P_E(y)^*x$$
and hence $P_E$ is an adjointable operator on $\OX$. Therefore $P_E$ defines (left
multiplication by) an element in the multiplier algebra $\M (\OX)$ \cite[Theorem~2.47]{RW}.  It is easy to check that $P_E^2 = P_E^* = P_E$ so that $P_E$ is a projection.  Furthermore, $P_E$ has properties (1), (2), and (3) in the statement of the lemma.

If we let $\{f_\lambda\}_{\lambda \in \Lambda}$ be an approximate unit for $A$, then a similar argument can be used to show $P_F(x) := \lim_\lambda \pi_X(f_\lambda, 0) x$ defines a projection in $\M ( \OX)$ with properties (4), (5), and (6).
\end{proof}

\begin{lemma} \label{EF-reg-X-reg}
If $E$ and $F$ are regular, then $R$ and $S$ are regular and $X$ is regular.
\end{lemma}

\begin{proof}
Since $E = R \otimes_B S$, we see that if $\phi_R(a) = 0$, then $$\phi_E(a) (r \otimes s) = \phi_R(a)(r) \otimes s = 0$$ for all $r \in R$ and $s \in S$.  Hence $\phi_E(a) = 0$ and by the injectivity of $\phi_E$ we have that $a=0$.  Thus $\phi_R$ is injective.  In addition, for any $a \in A$ we see that $\phi_E(a) \in \K(E)$. Since $\phi_E(a) = \phi_R(a) \otimes \textrm{Id}$ and $\phi_E$ is injective, Proposition~\ref{FMR-prop} implies that $\phi_R(a) \in \K(R)$.  Thus $R$ is regular.  Because $F = S \otimes_A R$, a similar argument shows that $S$ is regular.

Furthermore, since $X = S \oplus R$ and $\phi_X(a,b) = \phi_S(b) \oplus \phi_R(a)$ it is straightforward to show that $\phi_X$ is injective, and $\phi_X(a,b) \in \K(X)$ for all $(a,b) \in A \oplus B$.  Thus $X$ is regular.
\end{proof}

\begin{theorem} \label{SSE-correspondences}
Let $E$ be an essential, regular $C^*$-correspondence over a $C^*$-algebra $A$, and let $F$ be an essential, regular $C^*$-correspondence over a $C^*$-algebra $B$.  If $E$ is elementary strong shift equivalent to $F$, then $\mathcal{O}_E$ is Morita equivalent to $\mathcal{O}_F$.

In particular, if we write $E \cong R \otimes_B S$ and $F \cong S \otimes_A R$, with both $R$ and $S$ essential, then $\mathcal{O}_E$ and $\mathcal{O}_F$ are isomorphic to complementary full corners of $\OX$, where $X$ is the bipartite inflation of $S$ by $R$.
\end{theorem}

\begin{proof}
Since $E$ and $F$ are essential and elementary strong shift equivalent, we may use Corollary~\ref{choose-RS-essential} to write $E \cong R \otimes_B S$ and $F \cong S \otimes_A R$ with $R$ and $S$ essential.  Let $X$ be the bipartite dilation of $S$ by $R$ as defined in Definition~\ref{def-bipart-infl}, and let $(t_X, \pi_X) : (X, A \oplus B) \to \OX$ be a universal coisometric representation of $X$ into $\OX$.  Also let $(t, \pi) : (R \otimes_B S,A) \to \OX$ and $(t',\pi') : (S \otimes_A R,B) \to \OX$ be the coisometric representations defined in Definition~\ref{t-pi-X-def} and Definition~\ref{t-pi-prime-X-def}, respectively.  By Lemma~\ref{gauge-act-corner} and Lemma~\ref{t-pi-prime-X-def}, the representations $(t, \pi)$ and $(t', \pi')$ admit gauge actions.  Furthermore, $\pi$ and $\pi'$ are both injective since $\pi_X$ is injective.  Hence by the Gauge-Invariant Uniqueness Theorem these representations induce injective $*$-homomorphisms into $\OX$ and we have $\mathcal{O}_E \cong C^*(t,\pi)$ and $\mathcal{O}_F \cong C^*(t',\pi')$.

Let $P_E$ and $P_F$ be the projections in $\M (\OX)$ defined in Lemma~\ref{proj-mult-alg}.  We will prove that $P_E \OX P_E = C^*(t,\pi)$.  To begin, note that $$t (r \otimes s) = t_X (0,r) t_X(s,0) = P_E t_X (0,r) t_X(s,0) P_E  \in P_E \OX P_E$$ and since the elementary tensors span a dense subset of $R \otimes_B S$ we have that $\im t \subseteq P_E \OX P_E$.  Similarly, $\pi(a) = \pi_X(a,0) = P_E \pi(a,0) P_E \in P_E \OX P_E$ so that $\im \pi \subseteq P_E \OX P_E$.  Thus $C^*(t,\pi) \subseteq P_E \OX P_E$.

To see the reverse inclusion, note that $$\OX = \clspan \{ t_X(s_1, r_1) \ldots t_X(s_n, r_n) t_X(s'_m, r'_m)^* \ldots t_X(s'_1, r'_1)^* : m,n \geq 0 \}.$$  Thus to prove that $P_E \OX P_E \subseteq C^*(t,\pi)$ it suffices to show that $$P_E  t_X(s_1, r_1) \ldots t_X(s_n, r_n) t_X(s'_m, r'_m)^* \ldots t_X(s'_1, r'_1)^* P_E \in C^*(t,\pi).$$  To do this, we first notice the following equation holds:
\begin{align*}
& P_E t_X(s_1, r_1) t_X(s_2, r_2) \ldots t_X(s_n, r_n) \\
= & \ t_X(0, r_1) t_X(s_2, r_2) \ldots t_X(s_n, r_n) \\
= & \ t_X(0, r_1)P_F t_X(s_2, r_2) \ldots t_X(s_n, r_n) \\
= & \ t_X(0, r_1) t_X(s_2, 0) \ldots t_X(s_n, r_n) \\
= & \ t_X(0, r_1) t_X(s_2, 0) P_E \ldots t_X(s_n, r_n) \\
& \vdots \\
= & \ \begin{cases} t_X(0, r_1) t_X(s_2, 0) \ldots t_X(0, r_{n-1}) t_X(s_n, 0)P_E & \text{ if $n$ is even} \\
t_X(0, r_1) t_X(s_2, 0) \ldots t_X(s_{n-1}, 0) t_X(0, r_{n}) P_F & \text{ if $n$ is odd.}
 \end{cases}
\end{align*}
and then we consider three cases.

$ $

\noindent \textsc{Case 1:} $n$ and $m$ are both even

In this case 
\begin{align*}
& P_E  t_X(s_1, r_1) \ldots t_X(s_n, r_n) t_X(s'_m, r'_m)^* \ldots t_X(s'_1, r'_1)^* P_E \\
=& \  t_X(0, r_1) t_X(s_2, 0) \ldots t_X(0, r_{n-1}) t_X(s_n, 0) P_E P_E t_X(s'_m, 0)^* t_X(0,r'_{m-1})^* \\
& \qquad \ldots t_X(s'_2, 0)^* t_X(0, r'_1)^* \\
=& t_X(0, r_1) t_X(s_2, 0) \ldots t_X(0, r_{n-1}) t_X(s_n, 0) \\
& \qquad (t_X(0, r'_1) t_X(s'_2, 0) \ldots t_X(0, r'_{m-1}) t_X(s'_m, 0))^* \\
=& t (r_1 \otimes s_2) \ldots t( r_{n-1} \otimes s_n) (t (r'_1 \otimes s'_2) \ldots t( r'_{m-1} \otimes s'_m))^*
\end{align*}
 which is in $C^*(t,\pi)$.

$ $

\noindent \textsc{Case 2:} One of $m$ and $n$ is even and the other is odd.

First suppose that $n$ is odd and $m$ is even.  Then since $P_E$ and $P_F$ are orthogonal, we have that 
\begin{align*}
& P_E  t_X(s_1, r_1) \ldots t_X(s_n, r_n) t_X(s'_m, r'_m)^* \ldots t_X(s'_1, r'_1)^* P_E \\
=& \  t_X(0, r_1) t_X(s_2, 0) \ldots t_X(s_{n-1}) t_X(0, r_n) P_F P_E t_X(s'_m, 0)^* t_X(0,r'_{m-1})^* \\
& \qquad \ldots t_X(s'_2, 0)^* t_X(0, r'_1)^* \\
= & \ 0
\end{align*}
which is in $C^*(t,\pi)$.  The situation when $n$ is even and $m$ is odd is similar.

$ $

\noindent \textsc{Case 3:} $n$ and $m$ are both odd

Let $\{e_\lambda \}_{\lambda \in \Lambda}$ be an approximate unit for $A \oplus B$.  Then since $E$ and $F$ are regular, it follows from Lemma~\ref{EF-reg-X-reg} that $X$ is regular.  Hence we may write $\phi_X(e_\lambda) = \lim_j \sum_{k=1}^{N_j} \Theta^X_{(s^\lambda_{j,k},r^\lambda_{j,k}),(v^\lambda_{j,k},u^\lambda_{j,k})}$, and $$\pi_X(e_\lambda) = \psi_{t_X} (\phi_X(e_\lambda)) = \lim_j \sum_{k=1}^{N_j} t_X(s^\lambda_{j,k},r^\lambda_{j,k}) t_X(v^\lambda_{j,k},u^\lambda_{j,k})^*.$$  Then
\begin{align*}
& \ P_E t_X(s_1, r_1) \ldots t_X(s_n, r_n) t_X(s'_m, r'_m)^* \ldots t_X(s'_1, r'_1)^* P_E \\
=& \ \lim_\lambda P_E t_X(s_1, r_1) \ldots t_X(s_n, r_n) \pi(e_\lambda) t_X(s'_m, r'_m)^* \ldots t_X(s'_1, r'_1)^*P_E \\
=& \ \lim_\lambda \lim_j  \sum_{j=1}^{N_j} P_E t_X(s_1, r_1) \ldots t_X(s_n, r_n)  t_X(s^\lambda_{j,k},r^\lambda_{j,k})  \\
& \qquad \qquad \qquad \qquad t_X(v^\lambda_{j,k},u^\lambda_{j,k})^* t_X(s'_m, r'_m)^* \ldots t_X(s'_1, r'_1)^* P_E \\
\end{align*}
which is in $C^*(t,\pi)$ since it is a limit of sums of terms of the form  described in Case 1.

Thus we have shown that $C^*(t,\pi) = P_E \OX P_E$.  A similar argument shows that $C^*(t',\pi') = P_F \OX P_F$.  Thus $\mathcal{O}_E$ and $\mathcal{O}_F$ are isomorphic to the corners determined by $P_E$ and $P_F$, respectively.

To see that $C^*(t,\pi) = P_E \OX P_E$ is full, suppose that $\mathcal{I}$ is an ideal of $\OX$ containing $C^*(t,\pi)$.  Then $\mathcal{I}$ contains $\pi(a) = \pi_X(a,0)$ for all $a \in A$.  If $\{ f_\lambda \}_{\lambda \in \Lambda}$ is an approximate unit for $A$, then for any $s \in S$ we have that $\lim_\lambda (s,0) (e_\lambda, 0) = (s,0)$ so that $t_X(s,0) = \lim_\lambda t_X(s,0) \pi_X(e_\lambda,0) = \lim_\lambda t_X(s,0) \pi(e_\lambda)$ is in $\mathcal{I}$.  Furthermore, if $b \in B$, then since $S$ is regular by Lemma~\ref{EF-reg-X-reg}, we may write $\phi_S(b) = \lim_n \sum_{n=1}^{N_n} \Theta^S_{s_{n,k}, s'_{n,k}}$.  Because $X = S \oplus R$ as a right Hilbert $A \oplus B$-module, we see that $\phi_X(0,b) = \phi_S(b) \oplus 0 = \lim_n \sum_{n=1}^{N_n} \Theta^X_{(s_{n,k},0), (s'_{n,k},0)}$.  In addition, since $X$ is regular we may write $\pi(0,b) = \phi_X(0,b) = \lim_\lambda \sum_{n=1}^{N_n} t_X(s_{n,k}, 0) t_X(s'_{n,k}, 0)^*$, and thus $\pi(0,b)$ is in $\mathcal{I}$.  Hence for any $(a,b) \in A \oplus B$ we have that $\pi_X(a,b) = \pi_X(a,0) + \pi_X(0,b)$ is in $\mathcal{I}$.  But this implies that $\mathcal{I}$ is all of $\OX$.  Thus $C^*(t,\pi) = P_E \OX P_E$ is full.  A similar argument shows that $C^*(t',\pi') = P_F \OX P_F$ is full.

Finally, it follows from the relations in Lemma~\ref{proj-mult-alg} that $P_E + P_F = 1$ in $\M( \OX)$.  Thus the corners determined by $P_E$ and $P_F$ are complementary.  Since $\mathcal{O}_E$ and $\mathcal{O}_F$ are isomorphic to complementary full corners of $\OX$, it follows that $\mathcal{O}_E$ and $\mathcal{O}_F$ are Morita equivalent.
\end{proof}

\section{Non-essential $C^*$-correspondences} \label{non-ess-sec}

In Theorem~\ref{SSE-correspondences} we required that the $C^*$-correspondences $E$ and $F$ be essential.  It is unclear to the authors whether elementary strong shift equivalence of (not necessarily essential) regular $C^*$-correspondences $E$ and $F$ will always imply Morita equivalence of $\mathcal{O}_E$ and $\mathcal{O}_F$.  However, in this section we are able to prove that we may replace essentiality by the condition that the $C^*$-correspondence is over a unital $C^*$-algebra.

\begin{proposition} \label{J-X-J-ess}
Let $X$ be a $C^*$-correspondence over a $C^*$-algebra $A$.  If $A$ is unital, then $J_X = J_{X_\text{ess}}$.
\end{proposition}

\begin{proof}
Recall that $\phi_{X_\text{ess}}(a) = \phi_X(a)|_{X_\text{ess}}$ for all $a \in A$.  If $a \in A$, then for each $x \in X$ we see that $\phi_X(a)(x) = \phi_X(a) \phi_X(1)x = \phi_{X_\text{ess}}(a)(\phi_X(1)x)$.  It follows that $\ker \phi_X = \ker \phi_{X_\text{ess}}$.  

Additionally, if $a \in J(X)$, then $\phi_X(a) = \lim_n \sum_{k=1}^{N_n} \Theta^X_{x_{n,k}, y_{n,k}}$.  For any $z \in X_\text{ess}$ we may use the Hewitt-Cohen factorization Theorem \cite[Proposition~2.33]{RW} to write $z = \phi_X(b) w$ for $b \in A$ and $w \in X$.  We then have
\begin{align*}
\phi_X(a) z &= \phi_X(1) \phi_X(a) \phi_X(b) w \\
&= \phi_X(1) \lim_n \sum_{k=1}^{N_n} \Theta^X_{x_{n,k}, y_{n,k}} (\phi_X(b) w)\\
&= \phi_X(1) \lim_n \sum_{k=1}^{N_n} x_{n,k} \langle y_{n,k}, \phi_X(b) w \rangle_X \\
&= \lim_n \sum_{k=1}^{N_n} \phi_X(1) x_{n,k} \langle y_{n,k}, \phi_X(1^*)\phi_X(b) w \rangle_X \\
&= \lim_n \sum_{k=1}^{N_n} \phi_X(1) x_{n,k} \langle \phi_X(1)y_{n,k}, \phi_X(b) w \rangle_X \\
&= \lim_n \sum_{k=1}^{N_n} \Theta^{X_\text{ess}}_{\phi_X(1)x_{n,k}, \phi_X(1)y_{n,k}} (z)
\end{align*}
so that $\phi_{X_\text{ess}}(a) =  \lim_n \sum_{k=1}^{N_n} \Theta^{X_\text{ess}}_{\phi_X(1)x_{n,k}, \phi_X(1)y_{n,k}} \in \K(X_\text{ess})$ and $J(X) \subseteq J(X_\text{ess})$.  In addition, we see that if $a \in J(X_\text{ess})$, then $$\phi_{X_\text{ess}}(a) = \lim_n \sum_{k=1}^{N_n} \Theta^{X_\text{ess}}_{x'_{n,k}, y'_{n,k}}$$ and since $\phi_X(a)x = \phi_{X_\text{ess}}(a)(\phi_X(1)x)$ it is straightforward to show that $\phi_X(a) = \lim_n \sum_{k=1}^{N_n} \Theta^X_{x'_{n,k}, y'_{n,k}} \in \K(X)$.  Thus $J(X_\text{ess}) = J(X)$.  It follows that $J_{X_\text{ess}} = J_X$.
\end{proof}

\begin{proposition} \label{unital-ess-full-corn}
Let $X$ be a $C^*$-correspondence over a $C^*$-algebra $A$.  If $A$ is unital, then $\mathcal{O}_{X_\text{ess}}$ is isomorphic to a full corner of $\mathcal{O}_X$.  Consequently $\mathcal{O}_X$ is Morita equivalent to $\mathcal{O}_{X_\text{ess}}$.
\end{proposition}

\begin{proof}
Let $i : X_\text{ess} \hookrightarrow X$ be the inclusion map.  Let $(t_X, \pi_X) \to \mathcal{O}_X$ be the universal coisometric representation of $X$ into $\mathcal{O}_X$.  We define a representation $(t,\pi) : (X_\text{ess}, A) \to \mathcal{O}_X$ by setting $t := t_X \circ i$ and $\pi := \pi_X$.  It is straightforward to verify that $(t,\pi)$ is a representation.

To see that $(t, \pi)$ is coisometric, let $a \in J_{X_{\text{ess}}}$.  Then by Proposition~\ref{J-X-J-ess} we have that $a \in J_X$.  If we write $\phi_X(a) = \lim_n \sum_{k=1}^{N_n} \Theta^X_{x_{n,k}, y_{n,k}}$, then arguing as in Proposition~\ref{J-X-J-ess} shows that $\phi_{X_\text{ess}}(a) =  \lim_n \sum_{k=1}^{N_n} \Theta^{X_\text{ess}}_{\phi_X(1)x_{n,k}, \phi_X(1)y_{n,k}}$.  Thus 
\begin{align*}
\psi_t(\phi_{X_\text{ess}}(a)) &= \lim_n \sum_{k=1}^{N_n} t(\phi_X(1)x_{n,k}) t(\phi_X(1)y_{n,k})^* \\
&= \lim_n \sum_{k=1}^{N_n} \pi_X(1)t_X(x_{n,k}) t_X(y_{n,k})^*\pi_X(1) \\
&= \pi_X(1) \left( \lim_n \sum_{k=1}^{N_n}t_X(x_{n,k}) t_X(y_{n,k})^* \right) \pi_X(1) \\
&= \pi_X(1) \psi_{t_X}(\phi_X(a)) \pi_X(1) \\
&= \pi_X(1) \pi_X(a) \pi_X(1) \\
&= \pi_X(a) \\
&= \pi(a)
\end{align*}
so $(t, \pi)$ is coisometric.

It follows that $(t, \pi)$ induces a $*$-homomorphism $\rho_{(t,\pi)} : \mathcal{O}_{X_\text{ess}} \to \mathcal{O}_X$.  Since $\pi = \pi_X$ is injective, and since $(t, \pi)$ admits a gauge action (simply use the canonical gauge action of $\mathcal{O}_X$) it follows from the Gauge-Invariant Uniqueness Theorem that $\rho_{(t,\pi)}$ is injective and $\mathcal{O}_{X_\text{ess}} \cong C^*(t,\pi)$.

To see that $C^*(t,\pi)$ is a corner of $\mathcal{O}_X$ determined by the projection $\pi_X(1)$, simply note that 
\begin{align*}
&\pi_X(1) \Big( t_X(x_1) \ldots t_X(x_n) t_X(y_m)^* \ldots t_X(y_1)^* \Big) \pi_X(1) \\
& \qquad \qquad =  t_X(\phi_X(1)x_1) \ldots t_X(\phi_X(1)x_n) t_X(\phi_X(1)y_m)^* \ldots t_X(\phi_X(1)y_1)^* \\
& \qquad \qquad =  t(\phi_X(1)x_1) \ldots t(\phi_X(1)x_n) t(\phi_X(1)y_m)^* \ldots t(\phi_X(1)y_1)^*.
\end{align*}
Since the elements $t_X(x_1) \ldots t_X(x_n) t_X(y_m)^* \ldots t_X(y_1)^*$ span a dense subset of $\OX$ we see that $\pi_X(1) \OX \pi_X(1) \subseteq C^*(t,\pi)$.  

Furthermore, if $t(x_1) \ldots t(x_n) t(y_m)^* \ldots t(y_1)^* \in C^*(t,\pi)$, then because $\phi_X(1)x=x$ for $x \in X_\text{ess}$ we have that 
\begin{align*}
&t(x_1) \ldots t(x_n) t(y_m)^* \ldots t(y_1)^*  \\
& \qquad = \pi_X(1) t_X(x_1) \ldots t_X(x_n) t_X(y_m)^* \ldots t_X(y_1)^* \pi_X(1)
\end{align*}
which is in $\pi_X(1) \OX \pi_X(1)$.  Since the elements $t(x_1) \ldots t(x_n) t(y_m)^* \ldots t(y_1)^*$ span a dense subset of $C^*(t,\pi)$ we have that $C^*(t,\pi) \subseteq \pi_X(1) \OX \pi_X(1)$.  Thus $\pi_X(1) \OX \pi_X(1) = C^*(t,\pi)$.

Finally, to see that this corner is full, note that $\pi_X(1) \OX \pi_X(1)$ contains $\pi(A)$ and hence any ideal containing $\pi_X(1) \OX \pi_X(1)$ must be all of $\OX$.
\end{proof}

\begin{theorem}
Let $E$ be a regular $C^*$-correspondence over a $C^*$-algebra $A$, and let $F$ be an regular $C^*$-correspondence over a $C^*$-algebra $B$.  Suppose that either $E$ is essential or $A$ is unital.  Also suppose that either $F$ is essential or $B$ is unital.  If $E$ is elementary strong shift equivalent to $F$, then $\mathcal{O}_E$ is Morita equivalent to $\mathcal{O}_F$.
\end{theorem}

\begin{proof}
Since $E$ and $F$ are elementary strong shift equivalent we may write $E \cong R \otimes_A S$ and $F \cong S \otimes_B R$.  By Lemma~\ref{ess-across-otimes} we have that $E_\text{ess} \cong R_\text{ess} \otimes_A S_\text{ess}$ and $F_\text{ess} \cong S_\text{ess} \otimes_B R_\text{ess}$.   Since $E$ and $F$ are regular, it follows from Proposition~\ref{J-X-J-ess} that $E_\text{ess}$ and $F_\text{ess}$ are regular.  Therefore Theorem~\ref{SSE-correspondences} implies that $\mathcal{O}_{E_\text{ess}}$ is Morita equivalent to $\mathcal{O}_{F_\text{ess}}$.  Because $E$ is either essential or unital, we have that either $\mathcal{O}_E =  \mathcal{O}_{E_\text{ess}}$, or $\mathcal{O}_E$ is Morita equivalent to $\mathcal{O}_{E_\text{ess}}$ by Proposition~\ref{unital-ess-full-corn}.  Similarly, since $F$ is either essential or unital, we have that either $\mathcal{O}_F =  \mathcal{O}_{F_\text{ess}}$, or $\mathcal{O}_F$ is Morita equivalent to $\mathcal{O}_{F_\text{ess}}$ by Proposition~\ref{unital-ess-full-corn}.  Thus $\mathcal{O}_E$ is Morita equivalent to $\mathcal{O}_F$.
\end{proof}

\section{graph $C^*$-algebras: Examples and Counterexamples} \label{Examples and Counterexamples}

We use the conventions established in \cite{KPRR, KPR,BPRS, FLR, RS, BHRS} for graph $C^*$-algebras.  We also refer the reader to \cite{Rae} for a more comprehensive treatment of graph $C^*$-algebra theory --- although we warn the reader that the direction of the arrows in \cite{Rae} is ``opposite" of what is used in \cite{KPRR, KPR,BPRS, FLR, RS, BHRS} and of what is used here.

If $E = (E^0,E^1,r,s)$ is a graph, then the \emph{graph $C^*$-algebra} $C^*(E)$ is the universal $C^*$-algebra generated by a collection of mutually orthogonal projections $\{p_v : v \in E^0 \}$ together with a collection of partial isometries $\{ s_e : e \in E^1 \}$ with mutually orthogonal range projections that satisfy
\begin{enumerate}
\item $s_e^*s_e=p_{r(e)}$ \text{ for all $e \in E^1$}
\item $s_es_e^* \leq p_{s(e)}$ \text{ for all $e \in E^1$}
\item $p_v = \sum_{\{e : s(e)=v \} } s_es_e^*$ \text{ for all $v \in E^0$
with $0 < | s^{-1}(v) | < \infty$}.
\end{enumerate}

\noindent Alternatively, given a graph $E = (E^0,E^1,r,s)$ one may define a $C^*$-correspondence $X(E)$ over $A:= C_0(E^0)$ by letting 
\begin{equation*}
X(E) := \{ x : E^1 \to \mathbb{C} : \text{ the function } v \mapsto \sum_{
\{f \in E^1: r(f) = v \} } |x(f)|^2 \text{ is in $C_0(E^0)$} \ \}.
\end{equation*}
and giving $X(E)$  the operations 
\begin{align*}
(x \cdot a)(f) &:= x(f) a(r(f)) \text{ \quad for $f \in E^1$} \\
\langle x, y \rangle_{X(E)}(v) &:= \sum_{ \{ f \in E^1: r(f) = v \} }\overline{x(f)}y(f) \text{ \quad for $v \in E^0$} \\
(a \cdot x)(f) &:= a(s(f)) x(f) \text{ \quad for $f \in E^1$}.
\end{align*}
We call $X(E)$ the \emph{graph $C^*$-correspondence} associated to
$E$, and it is a fact that $\mathcal{O}_{X(E)} \cong C^*(E)$ \cite[Proposition~4.4]{FR}.  Thus the graph $C^*$-algebra may be thought of as the Cuntz-Pimsner algebra associated to the graph $C^*$-correspondence.  We refer the reader to \cite[\S 3]{MT} for a more detailed discussion and analysis of graph $C^*$-correspondences.

\subsection{Examples} \label{example-subsec}

We shall show how graph $C^*$-algebras give examples illustrating Theorem~\ref{SSE-correspondences}.

\begin{definition}
If $E = (E^0,E^1,r,s)$ is a graph, then the \emph{vertex matrix} of $E$ is the $E^0 \times E^0$ matrix $A_E$ with entries $$A_E(v,w) := \# \{ e \in E^0 : s(e) = v \text{ and } r(e) = w \}.$$
\end{definition}

Let  $E$ and $F$ be row-finite graphs with no sinks.  If the matrices $A_E$ and $A_F$ are elementary strong shift equivalent, then there are matrices $R$ and $S$ with non-negative entries for which $A_E = RS$ and $A_F = SR$.   It follows that $R$ must be a $E^0 \times F^0$ matrix, and $S$ must be a $F^0 \times E^0$ matrix.  In this case we may create a bipartite graph $G_{R,S}$ as follows:  We let $G_{R,S}^0 := E^0 \sqcup F^0$, and for $v \in E^0$ and $w \in F^0$ we draw $R(v,w)$ edges from $v$ to $w$, and $S(w,v)$ edges from $w$ to $v$.  (So, in particular, the vertex matrix of $G_{R,S}$ is $A_{G_{R,S}} = \left( \begin{smallmatrix} 0 & R \\ S & 0 \end{smallmatrix} \right)$ and $G_{R,S}$ is bipartite.)  It has been shown independently by Bates \cite[Theorem~5.2]{Bat} and by Drinen and Sieben \cite[Proposition~7.2]{DS} that $C^*(E)$ and $C^*(F)$ are isomorphic to complementary full corners of $C^*(G_{R,S})$, and thus Morita equivalent.

\begin{example} \label{graph-ex}
Let $E$ and $F$ be the following graphs.
$$\begin{matrix}
E & & \xymatrix{ v \ar[r]^b \ar@(dl,ul)^a & w \ar@(dr,ur)_c \\ } & & \qquad & &  F & &  \xymatrix{ x \ar[r]^e \ar@(dl,ul)^d & y \ar[r]^f &z \ar@(dr,ur)_g \\ }
\end{matrix}$$

Then we see that $A_E = \left( \begin{smallmatrix} 1 & 1 \\ 0 & 1 \end{smallmatrix} \right)$ and $A_F = \left( \begin{smallmatrix} 1 & 1 & 0 \\ 0 & 0 & 1 \\ 0 & 0 & 1 \end{smallmatrix} \right)$ are elementary strong shift equivalent by taking $R = \left( \begin{smallmatrix} 1 & 1 & 0 \\ 0 & 0 & 1 \end{smallmatrix} \right)$ and $S = \left( \begin{smallmatrix} 1 & 0 \\ 0 & 1 \\ 0 & 1 \end{smallmatrix} \right)$.  The bipartite graph $G_{R,S}$ is then equal to
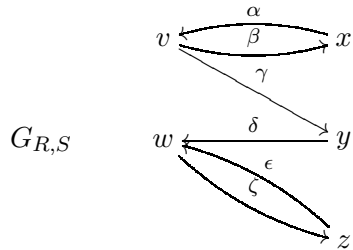
\begin{equation*}
\xymatrix{ & v \ar@/_/[rr]^\beta \ar[rrd]^\gamma & & x \ar@/_/[ll]_\alpha \\ G_{R,S} & w \ar@/_/[drr]^\zeta & & y \ar[ll]_\delta \\ & & & z \ar@/_/[llu]_\epsilon }
\end{equation*}
Also, $C^*(E)$ and $C^*(F)$ are isomorphic to complementary full corners of $C^*(G_{R,S})$.  In fact, if $\{ s_e, p_v\}$ is a generating Cuntz-Krieger $E$-family for $C^*(E)$ and if $\{S_e, P_v \}$ is a generating Cuntz-Krieger $G$-family for $C^*(G_{R,S})$, then the $*$-homomorphism that identifies $C^*(E)$ with a full corner of $C^*(G_{R,S})$ maps
\begin{equation*}
p_v \mapsto P_v, \quad p_w \mapsto P_w, \quad s_a \mapsto S_\beta S_\alpha, \quad s_b \mapsto S_\gamma S_\delta, \quad \text{ and } \quad s_c \mapsto S_\zeta S_\epsilon.
\end{equation*}
Also, if $\{t_e, q_v \}$ is a generating Cuntz-Krieger $F$-family for $C^*(F)$, then the $*$-homomorphism that identifies $C^*(F)$ with a full corner of $C^*(G_{R,S})$ maps
\begin{align*}
&q_x \mapsto P_x, \quad q_y \mapsto P_y, \quad q_z \mapsto P_z,\\ 
&t_d \mapsto S_\alpha S_\beta, \quad t_e \mapsto S_\alpha S_\gamma, \quad t_f \mapsto S_\delta S_\zeta, \quad \text{ and } \quad t_g \mapsto S_\epsilon S_\zeta.
\end{align*}
\end{example}

\begin{definition}
For a rectangular $I \times J$ matrix $R$ with non-negative entries we create a bipartite graph $G_R$ by defining $G_R^0 := I \sqcup J$ and for $i \in I$ and $j \in J$ we draw $R(i,j)$ edges from $i$ to $j$.

For this graph we may construct a $C^*$-correspondence $X_R$ from $A:= C_0(I)$ to $B := C_0(J)$ by setting
\begin{equation*}
X_R := \{ x : G_R^1 \to \mathbb{C} : \text{ the function } v \mapsto \sum_{
\{f \in G_R^1: r(f) = v \} } |x(f)|^2 \text{ is in $C_0(J)$} \ \}.
\end{equation*}
and giving $X_R$  the operations 
\begin{align*}
(x \cdot b)(f) &:= x(f) b(r(f)) \text{ \quad for $f \in G^1_R$} \\
\langle x, y \rangle_{X_R}(j) &:= \sum_{ \{ f \in G_R^1: r(f) = j \} }\overline{x(f)}y(f) \text{ \quad for $j \in J$} \\
(a \cdot x)(f) &:= a(s(f)) x(f) \text{ \quad for $f \in G_R^1$}.
\end{align*}
\end{definition}

\begin{example}
If $R$ and $S$ are the matrices in Example~\ref{graph-ex}, then $G_R$ and $G_S$ are the following graphs:
$$
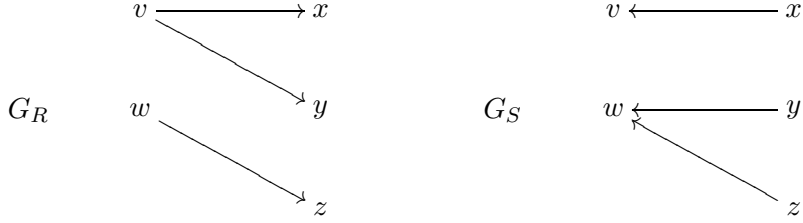
\begin{matrix}
 \xymatrix{ & v \ar[rr] \ar[rrd] & & x \\ G_R & w \ar[drr] & & y \\ & & & z }
 & \qquad & & \xymatrix{ & v & & x \ar[ll] \\ G_S & w & & y \ar[ll] \\ & & & z \ar[llu] }
\end{matrix}$$
\end{example}

This relates to Theorem~\ref{SSE-correspondences} in the following way.  Since $E$ and $F$ are row-finite with no sinks, it follows that $X(E)$ and $X(F)$ are regular \cite[Remark~3.3]{MT}.  In addition, graph $C^*$-correspondences are always essential \cite[\S 3]{MT}.

The fact that $A_E = RS$ and $A_F = SR$ implies that $X(E) \cong X(G_R) \otimes_B X(G_S)$ and $X(F) \cong X(G_S) \otimes_A X(G_R)$.  Hence the graph $C^*$-correspondences $X(E)$ and $X(F)$ are elementary strong shift equivalent.  The bipartite inflation of $X(G_S)$ by $X(G_R)$ is equal to $X(G_{R,S})$.  Thus Theorem~\ref{SSE-correspondences} implies that $\mathcal{O}_{X(E)} \cong C^*(E)$ and $\mathcal{O}_{X(F)} \cong C^*(F)$ are isomorphic to complementary full corners of $\mathcal{O}_{X(G_{R,S})} \cong C^*(G_{R,S})$.  In this way we recover \cite[Theorem~5.2]{Bat} and \cite[Proposition~7.2]{DS} as special cases of Theorem~\ref{SSE-correspondences}.

\subsection{Counterexamples}
We shall use graph $C^*$-algebras to provide counterexamples to generalizations of the statement of Theorem~\ref{SSE-correspondences}.

We have already mentioned in \S \ref{non-ess-sec} that the authors are unsure whether Theorem~\ref{SSE-correspondences} remains true if $E$ and $F$ are not essential.  On the other hand, the condition that $E$ and $F$ are regular is necessary.  When we impose the condition that a $C^*$-correspondence be regular, we require the left action to be injective as well as act by compact operators.  We shall show that both of these conditions are necessary in Theorem~\ref{SSE-correspondences}.

If we let $E_1$ and $E_2$ be the following graphs
$$\begin{matrix}
E_1 & & \xymatrix{ v & w \ar[l] \ar[r]  & x \\ } & & \qquad & &  E_2 & &  \xymatrix{ y \ar[r] & z \\ }
\end{matrix}$$
then $A_{E_1} = \left( \begin{smallmatrix} 0 & 0 & 0 \\ 1 & 0 & 1 \\ 0 & 0 & 0\end{smallmatrix} \right)$ and $A_{E_2} = \left( \begin{smallmatrix} 0 & 1 \\ 0 & 0 \end{smallmatrix} \right)$ are elementary strong shift equivalent by taking $R = \left( \begin{smallmatrix} 0 & 0 \\ 0 & 1 \\ 0 & 0 \end{smallmatrix} \right)$ and $S = \left( \begin{smallmatrix} 0 & 1 & 0 \\ 1 & 0 & 1 \end{smallmatrix} \right)$.  Reasoning as in \S \ref{example-subsec} shows that the graph $C^*$-correspondences $X(E_1)$ and $X(E_2)$ are elementary strong shift equivalent.  Furthermore, $X(E_1)$ and $X(E_2)$ are essential (as are all graph $C^*$-correspondences).  In addition, since $E_1$ and $E_2$ are row-finite, their left actions act as compact operators, and we have $J(X(E_1)) = C_0(E_1^0)$ and $J(X(E_2)) = C_0(E_2^0)$.  However, since each of $E_1$ and $E_2$ has sinks, the left actions of the associated graph $C^*$-correspondences are not injective.  Thus neither $X(E_1)$ nor $X(E_2)$ is regular.  In addition, $\mathcal{O}_{X(E_1)} \cong C^*(E_1)$ is not Morita equivalent to $\mathcal{O}_{X(E_2)} \cong C^*(E_2)$ because $C^*(E_1)$ contains two proper ideals and $C^*(E_2)$ is simple.  This shows that we cannot remove the regularity condition in Theorem~\ref{SSE-correspondences}.

Moreover, there is an example, described in \cite[Example~5.4]{Bat}, which shows that there exist non-row-finite graphs $E_1$ and $E_2$ with no sinks that have graph $C^*$-correspondences that are elementary strong shift equivalent but have $C^*$-algebras that are not Morita equivalent.  Thus one needs the left action to be both injective and act as compact operators.

Finally, we mention some natural questions that arise when one considers elementary strong shift equivalence of $C^*$-correspondences.  We have seen that elementary strong shift equivalence of $C^*$-correspondences implies Morita equivalence of the associated Cuntz-Pimsner algebras.  It is natural to ask whether this equivalence holds at higher levels --- in particular, at the level of Toeplitz algebras, or at the level of $C^*$-correspondences.  Thus there are three natural questions one can ask.

\

\noindent Let $E$ and $F$ be essential, regular $C^*$-correspondences.

\smallskip

\smallskip

\noindent \textsc{Question 1:}  If $E$ and $F$ are elementary strong shift equivalent, then is it necessarily the case that $E$ and $F$ are Morita equivalent (as defined in \cite{MS2})?

\smallskip

\noindent \textsc{Question 2:}  If $E$ and $F$ are elementary strong shift equivalent, then is it necessarily the case that the Toeplitz algebras $\mathcal{T}_E$ and $\mathcal{T}_F$ are Morita equivalent?

\smallskip

\noindent \textsc{Question 3:}  If $E$ and $F$ are elementary strong shift equivalent, then is it necessarily the case that the Cuntz-Pimsner algebras $\mathcal{O}_E$ and $\mathcal{O}_F$ are Morita equivalent?

\smallskip

\smallskip

We have seen that Theorem~\ref{SSE-correspondences} provides an affirmative answer to Question~3.  In addition, the questions asked above are successively weaker in the following sense:  If $E$ and $F$ are Morita Equivalent as $C^*$-correspondences, then it follows that $\mathcal{T}_E$ and $\mathcal{T}_F$ are Morita Equivalent.  Furthermore, since the Cuntz-Pimsner algebra is a quotient of the Toeplitz algebra by a certain ideal, we see that if the Morita equivalence between $\mathcal{O}_E$ and $\mathcal{O}_F$ takes the appropriate ideal to the appropriate ideal, then $\mathcal{O}_E$ and $\mathcal{O}_F$.

When the authors began this project, they intended to prove a theorem that would provide an affirmative answer to Question~1, and then obtain affirmative answers to Question~2 and Question~3 as corollaries by using the arguments of the previous paragraph.  However, upon deeper investigation it appears that Question~2 and Question~3 have negative answers.  In particular, the Morita equivalence can only be guaranteed to hold at the level of Cuntz-Pimsner algebras, and not at the level of Toeplitz algebras or $C^*$-correspondences.

To see this, let $E$ and $F$ be the following graphs.  

\

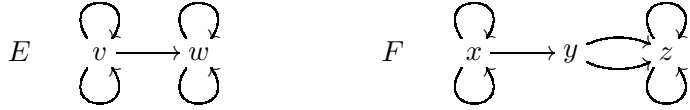
$$\begin{matrix}
E & & \xymatrix{ v \ar[r] \ar@(ul,ur) \ar@(dl,dr) & w \ar@(ul,ur) \ar@(dl,dr) \\ } & & \qquad & &  F & &  \xymatrix{ x \ar[r] \ar@(ul,ur) \ar@(dl,dr) & y \ar@/^/[r] \ar@/_/[r] &z \ar@(ul,ur) \ar@(dl,dr) \\ } \\
& & & & & & & & 
\end{matrix}$$

Then we see that $A_E = \left( \begin{smallmatrix} 2 & 1 \\ 0 & 2 \end{smallmatrix} \right)$ and $A_F = \left( \begin{smallmatrix} 2 & 1 & 0 \\ 0 & 0 & 2 \\ 0 & 0 & 2 \end{smallmatrix} \right)$ are elementary strong shift equivalent by taking $R = \left( \begin{smallmatrix} 2 & 1 & 0 \\ 0 & 0 & 2 \end{smallmatrix} \right)$ and $S = \left( \begin{smallmatrix} 1 & 0 \\ 0 & 1 \\ 0 & 1 \end{smallmatrix} \right)$.  Reasoning as in \S \ref{example-subsec} shows that the graph $C^*$-correspondences $X(E)$ and $X(F)$ are elementary strong shift equivalent.  Also, $X(E)$ and $X(F)$ are essential, and since $E$ and $F$ are row-finite with no sinks, it follows that $X(E)$ and $X(F)$ are regular.

As discussed in \cite[\S 3]{MT} and \cite[Theorem~3.7]{MT} the Toeplitz algebra of $X(E)$ is the $C^*$-algebra of the graph formed by outsplitting $E$ at all of its vertices.  Similarly for $F$.  Thus if we let $\widetilde{E}$ and $\widetilde{F}$ be the following graphs

\

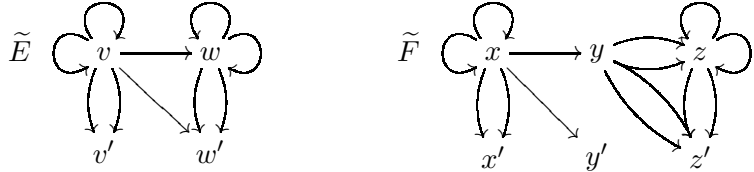
$$\begin{matrix}
\widetilde{E} & & \xymatrix{ v \ar[r] \ar@(ul,ur) \ar@(ul,dl) \ar@/^/[d] \ar@/_/[d] \ar[dr] & w \ar@(ur,ul) \ar@(ur,dr) \ar@/_/[d] \ar@/^/[d] \\  v' & w' } & & \qquad & &  \widetilde{F} & &  \xymatrix{ x \ar[r] \ar@(ul,ur) \ar@(ul,dl) \ar@/^/[d] \ar@/_/[d] \ar[dr] & y  \ar@/^/[r] \ar@/_/[r]  \ar@/^/[rd] \ar@/_/[rd] &z \ar@(ur,ul) \ar@(ur,dr)  \ar@/^/[d] \ar@/_/[d] \\ x' & y' & z'}
\end{matrix}$$
then $\mathcal {T}_{X(E)} \cong C^*(\widetilde{E})$ and $\mathcal{T}_{X(F)} \cong C^*(\widetilde{F})$.  However, since the proper saturated, hereditary subsets of $\widetilde{E}$ are $$ \emptyset \qquad \{ w' \} \qquad \{w, w' \} \qquad \{ v' \} \qquad \{v', w' \} \qquad \{v', w, w' \}$$
and because $\widetilde{E}$ satisfies Condition~(K) we see that $C^*(\widetilde{E})$ has exactly 6 proper ideals.  In addition, since the saturated, hereditary subsets of $\widetilde{F}$ are 
\begin{align*}
\emptyset \quad &\{ z' \} \quad \{y', z' \} \quad \{y, z, z' \} \quad \{y, y', z, z' \} \quad \{ x' \} \quad \{ x', z' \} \\
&\quad \{x', y', z' \} \quad \{x', y, z, z' \} \quad \{x', y, y', z, z' \}
\end{align*}
and because $\widetilde{F}$ satisfies Condition~(K) we see that $C^*(\widetilde{F})$ has exactly 10 proper ideals.  Thus $\mathcal {T}_{X(E)} \cong C^*(\widetilde{E})$ and $\mathcal{T}_{X(F)} \cong C^*(\widetilde{F})$ are not Morita equivalent, and this provides a negative answer to Question 2.  Moreover, it provides a negative answer to Question 1, since Morita equivalence of $C^*$-correspondences implies Morita equivalence of the associated Toeplitz algebras.

\begin{remark}
We conclude with a thought which motivated us at the outset, but which we could not verify.  Suppose that $E$ is an essential, regular $C^*$-correspondence over a $C^*$-algebra $A$. Then for every $n \in \mathbb{N}$, the map $T \mapsto T\otimes \textrm{Id}_E$ that embeds $\mathcal{L}(E^{\otimes n})$ in $\mathcal{L}(E^{\otimes (n+1)})$ carries $\mathcal{K}(E^{\otimes n})$ into $\mathcal{K}(E^{\otimes (n+1)})$.  Let $\mathfrak{A}$ denote the inductive limit $\dirlim \mathcal{K}(E^{\otimes n})$ and let $\mathcal{E}: = E \otimes_A \mathfrak{A}$. Then $\mathcal{E}$ is an invertible correspondence over $\mathfrak{A}$ in the sense that $\mathcal{E}$ is an imprimitivity bimodule (from $\mathfrak{A}$ to $\mathfrak{A}$) and, as is shown in \cite[Theorem 2.5]{Pim}, the Cuntz-Pimsner algebra $\mathcal{O}_{\mathcal{E}}$ is isomorphic to $\mathcal{O}_E$.  Suppose too that $F$ is an essential, regular $C^*$-correspondence over the $C^*$-algebra $B$ and let $\mathfrak{B}$ and $\mathcal{F}$ be the analogous inductive limit and invertible correspondence.  One of our initial approaches to proving Theorem \ref{SSE-correspondences} was to try to prove that if $E$ and $F$ are strong shift equivalent then $\mathcal{E}$ and $\mathcal{F}$ are Morita equivalent in the sense of \cite{MS2}. While the implication still seems plausible, we are unable to decide whether it is true or false.  It seems like the ``right'' conjecture to make in view of Williams's theorems.  In fact, one is enticed to speculate on its converse, too: If $\mathcal{E}$ and $\mathcal{F}$ are Morita equivalent in the sense of \cite{MS2}, then are $E$ and $F$ strong shift equivalent? An ``if and only if'' theorem would indeed be a perfect analogue of Williams's theorems.
\end{remark}

\end{document}